\documentclass[a4paper,12pt]{amsart}
\usepackage[utf8]{inputenc}
\usepackage[T1]{fontenc}
\usepackage[UKenglish]{babel}
\usepackage[a4paper,margin=28mm]{geometry}
\usepackage{verbatim}
\allowdisplaybreaks[4]
\usepackage{times}
\usepackage{dsfont,mathrsfs}
\usepackage{latexsym}
\usepackage{subfigure}
\usepackage{mathtools}
\usepackage{amsmath,amsfonts,amssymb,amssymb,amsfonts,amsthm,color}
\usepackage{mathrsfs}
\usepackage{graphicx}  
\usepackage{tikz, subfigure}  
\usepackage{url}
\usepackage{enumitem}
\usepackage{float}
\usepackage{hyperref} 
\hypersetup{hypertex=true,
	colorlinks=true,
	linkcolor=blue,
	anchorcolor=blue,
	citecolor=blue,
	pdftitle={FCLT EMscheme},
	pdfpagemode=}
\usepackage{xcolor}
\usepackage[author-year]{amsrefs}

\theoremstyle{plain}
\newtheorem{theorem}{Theorem}[section]

\newtheorem{lemma}[theorem]{Lemma}
\newtheorem{assumption}[theorem]{Assumption}

\newtheorem{problem}[theorem]{Problem}
\newtheorem{remark}{Remark}[section]

\theoremstyle{definition}

\numberwithin{equation}{section}
\renewcommand\labelenumi{\textup{\alph{enumi})}}
\renewcommand\theenumi\labelenumi

\usepackage{marginnote}

\begin{document}
	\title[FCLT for EM scheme with decreasing step sizes]{Functional central limit theorem for Euler--Maruyama scheme with decreasing step sizes}
\date{\today}
\author[Q.~Pei]{Qiyang Pei}
\address[Q.~Pei]{Department of Mathematics, Faculty of Science, University of Macau, Macau S.A.R., China.}
\email{yc37431@um.edu.mo, qiyangpei@gmail.com}
\author[L.~Xu]{Lihu Xu}
\address[L.~Xu]{Department of Mathematics, Faculty of Science, University of Macau, Macau S.A.R., China.}
\email{lihuxu@um.edu.mo}

\maketitle
\pagestyle{headings}
\begin{abstract}
We consider the Euler--Maruyama (EM) scheme of a family of dissipative stochastic differential equations (SDEs), whose step sizes $\eta_{1}\ge\eta_{2}\ge\eta_{3}\ge \cdots$ are decreasing, and prove that the EM scheme weakly converges to a time-changed Brownian motion (BM) $\{B_{a(t)}^{1}\}_{0\le t\le 1}$ rather than the standard BM $\{B_{t}^{1}\}_{0\le t\le1}$, where the increasing function $a(t)$  depends on the choice of $\{\eta_{k}\}_{k\ge1}$, for instance, $a(t)=t^{1+\beta}$ if $\eta_{k}=k^{-\beta}$.

Compared to the EM scheme with constant step size, there are substantial differences as follows:
\begin{itemize}
	\item[(i)] the EM time series is inhomogeneous and weakly converges to the  ergodic measure at a polynomial rate; 
	\item[(ii)] we have a quantity $T_{n}=\frac{1}{\eta_{1}}+\cdots+\frac{1}{\eta_{n}}$ which roughly measures the dependence of the EM time series;
	\item[(iii)] the normalized factor in the CLT is $n/\sqrt{T_{n}}$ rather than $\sqrt{n}$, in particular, $n/\sqrt{T_{n}} \asymp n^{(1-\beta)/2}$ when $\eta_{k}=k^{-\beta}$ with $\beta\in(0,1)$;  
	\item[(iv)] in the critical choice $\eta_{k}=1/k$, we have $n/\sqrt{T_{n}}=O(1)$ and thus conjecture that the CLT and FCLT do not hold. This conjecture has been verified by simulations. 
\end{itemize}
A key distinction arises between the constant and decreasing step size implementations of the EM scheme. Under a constant step size, the time series is time-homogeneous, which allows us to use a stationary initialization. This trick automatically eliminates several complex terms in the proof. Conversely, the time series generated by the EM scheme with decreasing step sizes forms an inhomogeneous Markov chain. To manage the analogous difficult terms in this case, especially when the test function $h$ is Lipschitz, we must instead establish a bound for the Wasserstein--2 distance. This technique for handling the inhomogeneous case could be of independent interest. 
\end{abstract} 

\section{Introduction}\label{SectionIntro} 
We consider the following stochastic differential equation (SDE),
\begin{align}
	\label{def SDE}
	{\rm d}X_{t}=b(X_{t}){\rm d}t+  \sigma  {\rm d}B_{t},\ X_{0}=x \in\mathbb{R}^{d},
\end{align}
where $\{B_{t}\}_{t\ge0}$ is $d$-dimensional standard Brownian motion, $\sigma$ is an invertible $d \times d$ matrix, and the drift $b:\mathbb{R}^{d}\to\mathbb{R}^{d}$ above is assumed to satisfy Assumption \ref{Asp Coeff}, which ensures that SDE \eqref{def SDE} admits a unique strong solution and a unique invariant measure $\pi$; see more details in \ocites{Mao08,LTX22}.

In practice, the invariant measure $\pi$ of SDE \eqref{def SDE} is simulated by unadjusted Langevin algorithm (ULA), using the Euler--Maruyama (EM) scheme of \eqref{def SDE}. The EM scheme of SDE \eqref{def SDE}  reads as
\begin{align} \label{def EMDis}\theta_{k+1}=\theta_{k}+\eta_{k+1}b(\theta_{k})+\sqrt{\eta_{k+1}}\sigma\xi_{k+1},\ \theta_{0}=x ,\quad \xi_{k+1}\stackrel{\text{i.i.d.}}{\sim}\mathcal{N}(0,I_{d}), 
\end{align}
where $\mathcal{N}(0,I_{d})$ denotes the $d$-dimensional standard normal distribution, and $\{\eta_{k}\}_{k \ge 1}$ is a decreasing sequence satisfying assumptions in Section \ref{SectionMainResults}. 

It is easy to see that $\{\theta_{k}\}_{k\ge0}$ is an time-inhomogeneous Markov chain for non-constant $\{\eta_{k}\}_{k\ge0}$. Define  
\begin{align}
	\label{def t_{n}}t_{0}\coloneqq0,\quad \text{ and }\quad t_{k}\coloneqq\eta_{1}+\cdots+\eta_{k},\quad k\ge 1,
\end{align}
and define the empirical measure associated with $\{\theta_{k}\}_{k\ge0}$ as the following: 
\begin{align} \label{def Pi_n}
	\Pi_{n}\coloneqq\frac{1}{n}\sum_{k=0}^{n-1} \delta_{\theta_{k}},\quad n \ge 1,
\end{align}
where $\delta_{\theta_{k}}$ is the delta measure at the point $\theta_{k}$. It has been shown in \ocite{LWX25}*{Theorem 2.1} that 
the distribution of $\{\theta_{k}\}_{k\ge 0}$ converges to $\pi$ and that the empirical measure $\Pi_n$ weakly converges to $\pi$.  

We shall prove a central limit theorem (CLT) and a functional central limit theorem (FCLT) for the EM scheme sequence $\{\theta_{k}\}_{k\ge 0}$, more precisely, as $n\to\infty$, determining weak limits of
\begin{align}
	\label{def prelim sum} 
	  \frac{n}{\sqrt{T_{n}}}[\Pi_{n}(h)-\pi(h)]
\end{align}
and
\begin{align} \label{def prelim path}
  \left\{\frac{\lfloor n t\rfloor }{\sqrt{T_{n}}}[\Pi_{\lfloor nt \rfloor}(h)-\pi(h)]\right\}_{t \in [0,1]}  \quad {\rm in}\ D[0,1],
\end{align}
where 
\begin{align} \label{def T_n}
	T_{n}\coloneqq\frac{1}{\eta_{1}}+\cdots+\frac{1}{\eta_{n}},
\end{align}
and $h$ is a Lipschitz function, $\lfloor x\rfloor$ denotes the integer part of $x \in \mathbb{R}$, and $D[0,1]$ denotes the c\`{a}dl\`{a}g function space from $[0,1]$ to $\mathbb{R}$.

\subsection{Related work}

\ocite{LP02} proved CLTs for weighted empirical average  $\sqrt{t_{n}}\nu_{n}^{\eta}(f)\coloneqq\frac{1}{\sqrt{t_{n}}}\sum_{k=1}^{n}\eta_{k}f(\theta_{k-1})$ for EM schemes $\{\theta_{k}\}_{k\ge 1}$ with test functions $f=b$ and $f=\mathcal{A}\varphi$, where $\mathcal{A}$ denotes the generator of $\eqref{def SDE}$.   \ocite{LTX22}  established a CLT with $h \in C_{b}^2(\mathbb {R}^{d};\mathbb{R})$ for the EM scheme of \eqref{def SDE} whose step sizes are constant; this result is extended  to the case with the test function $h \in\text{{\rm Lip}}(\mathbb{R}^{d};\mathbb{R})$ by \ocite{DFL26}. 
\ocite{LR23} established an FCLT for stochastic gradient Langevin dynamics (SGLD) with constant step size, which converges to a standard Brownian motion. \ocite{YBVE21} derived a CLT for stochastic gradient descent (SGD) with constant step size by applying \ocite{MT09}*{Theorem 17.0.1}. \ocite{LXY24} established CLTs for the average of the outputs of SGD and momentum SGD with decreasing step sizes. Related results on moderate deviations can be found in \ocites{DFL26,FHX24,GJW18,QY22}.  \ocite{Jin25} later analyzed the backward EM scheme, obtained by replacing $b(\theta_{k})$ in \eqref{def EMDis} with $b(\theta_{k+1})$, also under constant step size, and proved a CLT.

When $b(x)=-\nabla U(x)$ with $U$ being a potential function and $\sigma=\sqrt{2} I_{d}$ with $I_d$ the $d \times d$ identity matrix, $\pi$ admits a probability density proportional to $e^{-U(x)}$, and \eqref{def EMDis} reduces to the corresponding ULA. See, for example \ocites{Dal17,DM19}, for the applications of ULA. When $U$ is convex outside a ball and $\eta_{k}=\eta$ for all $k$, \ocite{MMS20} proved that ULA converges to $\pi$ at a speed $\eta^{1/4}$. This result has been extended to tamed ULA in \ocites{LNSZ24,NNZ25}.

\subsection{Contributions and approach} 
Our  main results are a CLT and an FCLT. Different from the classical CLT, our scaling number is $n/\sqrt{T_{n}}$ rather than $\sqrt n$. This is because the time series $\{\theta_{k}\}_{k \ge 0}$
is an inhomogeneous Markov chain and more and more correlated as $k \rightarrow \infty$. To see this, let us take $\eta_{k}=k^{-\beta}$ with $\beta \in [0,1]$. We have\begin{align*}
	\frac{n}{\sqrt{T_{n}}}\asymp n^{(1-\beta)/2}.
\end{align*} When $\beta=0$, $n/\sqrt{T_{n}} \asymp n^{1/2}$ which is the order of time-homogeneous time series. Because the correlation of $\{\theta_{k}\}_{k \ge 0}$ increases with respect to $\beta$,  the scaling number $n/\sqrt{T_{n}}$ of CLT decreases. When $\beta=1$, we arrive at a critical point and conjecture that the CLT does not hold true. This has been verified by simulation below. To prove the CLT, we use the decomposition \eqref{def Decomposition} by solving the Poisson equation \eqref{def Poisson equation}, so that $[\Pi_{n}(h)-\pi(h)]$ is written as sum of a martingale and some negligible remainders.

As for the FCLT, unlike the classical case, the limit of our FCLT is a Brownian motion with a time change $a(t)$, i.e. $\{B_{a(t)}\}_{t \in [0,1]}$, where 
\begin{align*}
	a :[0,1] \to \mathbb{R}_{+}.
\end{align*}
The form of this $a(t)$ depends on the choice of step sizes $\{\eta_{k}\}_{k \ge 1}$. 
\ocite{McL74}*{Theorem 3.2} gives the classical criterion for FCLT, which cannot be directly applied. Following the same argument of proving \ocite{McL74}*{Theorem 3.2} with a minor modification, we  generalized this criterion, which yields our FCLT immediately. 

Since the Markov chain $\{\theta_{k}\}_{k \ge 0}$ of \eqref{def EMDis} is inhomogeneous, it converges to the ergodic measure polynomially rather than exponentially; see \ocite{LWX25}. This makes the techniques based on the homogeneous Markov chain such as the stationary initialization in \ocites{FSX19,LTX22} or techniques based on $\alpha$-mixing in \ocite{LR23} no longer apply.  More precisely, thanks to the homogeneous Markov chain, \ocites{LTX22,DFL26} applied a stationary initialization trick so that 
$\{\theta_{k}\}_{k \ge 0}$ is stationary, making several terms in the decomposition \eqref{def Decomposition} below automatically vanish, and that bounding the Wasserstein-1 error $\mathbf{W}_{1}(\theta_{k},X_{t_{k}})$ is sufficient. However, in our case, it is difficult to handle these terms and we need to bound the error between $\theta_{k}$ and $X_{t_{k}}$ in the Wasserstein--2 ($\mathbf{W}_{2}$) distance. By the coupling established in \ocite{Ebe16}, we show that   $\mathbf{W}_{2}(\theta_{k}, X_{t_{k}})$ is bounded by $\eta_{k}^{1/4}$, which is consistent with the result in \ocite{MMS20}.

The second difficulty in our case arises from the Lipschitz test function family, which is much larger than  $C_{b}^{2}(\mathbb{R}^{d};\mathbb{R})$. The corresponding Poisson equation has a worse regularity. This makes bounding the error terms in the decomposition \eqref{def Decomposition} much more involved.

\subsection{Notations}
For two real numbers $a,b$, we denote $a\lor b\coloneqq\max\{a,b\}$ and for $a>0$ define $\lfloor a\rfloor\coloneqq\sup\{n\in\mathbb{Z}:n\le a\}$. For a finite set $A$, we write $\sharp A$ to denote its cardinality. For sums and products of sequence $\{a_{i}\}_{i\ge k}$, by convention, we set $\sum_{i=k}^{k-1}a_{i}=0$ and $\prod_{i=k}^{k-1}a_{i}=1$.

Let $I_{d}\in\mathbb{R}^{d\times d}$ be $(d\times d)$-identity matrix. For the symmetric matrices $A,B\in\mathbb{R}^{d\times d}$, write $A\preceq B$ if $B-A$ is positive semi-definite. The notation $\vert\cdot\vert$ denotes the Euclidean norm for $\mathbb{R}^{d}$-vectors and  the Frobenius norm for matrices. We write $f(n)\asymp g(n)$ if there are constants $C_{1},C_{2}>0$ such that for all $n\ge n_{0}$, $C_{1}f(n)\le g(n)\le C_{2}f(n)$. Sometimes we use $f(n)\lesssim g(n)$ if $f(n)\le Cg(n)$ for some constant $C>0$ independent of $n$.

For a function $f:\mathbb{R}^{d}\to\mathbb{R}$ that is sufficiently smooth, we denote its gradient and Hessian at $x$ by $\nabla f(x)\in\mathbb{R}^{d}$ and $\nabla^{2}f(x)\in\mathbb{R}^{d\times d}$ respectively. We denote by $C_{b}^{k}(\mathbb{R}^{d};\mathbb{R})$ the space of the bounded continuous functions from $\mathbb{R}^d$ to $\mathbb{R}$ whose $1$-st, ... ,$k$-th order derivatives are all bounded continuous. Denote by $\textbf{\rm{Lip}}(\mathbb{R}^{d};\mathbb{R})$ the collection of Lipschitz functions from $\mathbb{R}^{d}$ to $\mathbb{R}$. For a bounded continuous function $f:\mathbb{R}^{d}\to\mathbb{R}$, we define the uniform norm $\Vert f\Vert_{\infty}\coloneqq\sup_{x\in\mathbb{R}^{d}}\vert f(x)\vert$.

If a random variable $X$ has distribution (or law) $\mu$, we write $X\sim \mu$ or $\text{\rm{Law}}(X)=\mu$; if two random variables $X$ and $Y$ have the same law, we write $X\stackrel{d}{=}Y$. The notation $X_{n}\Rightarrow X$ denotes that  random variable sequence $\{X_{n}\}_{n}$ converges weakly to  $X$. We denote by $\mathcal{N}(\mu,\Sigma)$ the Gaussian distribution with mean $\mu$ and covariance $\Sigma$. For a measure $\pi$ and a measurable function $f$ defined on $\mathbb{R}^{d}$, we write 
\begin{align*}
	\pi(f)=\int_{\mathbb{R}^{d}} f(x) \pi({\rm d}x).
\end{align*}
The conditional expectation $\mathbb{E}[\cdot\vert \theta_{k}]$ is abbreviated as $\mathbb{E}_{k}[\cdot]$. Denote by $D[0,T]$ the Skorohod $D$--space on $[0,T]$ for $T>0$.

The Wasserstein-$p$ distance between probability measures $\mu$ and $\nu$ is defined as \begin{align*}
	\textbf{\rm{W}}_{p}(\mu,\nu)=\inf_{\gamma\in\textbf{C}(\mu,\nu)}[\mathbb{E}_{(X,Y)\sim\gamma}\vert X-Y\vert^{p}]^{1/p},
\end{align*} 
where $\textbf{C}(\mu,\nu)$ denotes the set of all couplings with marginal distributions $\mu$ and $\nu$. The Kantorovich distance of two distributions $\mu$ and $\nu$ on a metric space $(S,\rho)$ is defined by \begin{align}\label{defWrho}
	\mathbf{W}_{\rho}(\nu,\mu)=\inf_{\gamma\in \textbf{C}(\nu,\mu)}\int\rho(x,y)\gamma({\rm d}x,{\rm d}y).
\end{align}
For random variables  $X\sim \mu$ and $Y\sim\nu$, we will also use the notations $W_{p}(X,Y)\coloneqq W_{p}(\mu,\nu)$ and $W_{\rho}(X,Y)\coloneqq W_{\rho}(\mu,\nu)$.

To study EM scheme \eqref{def EMDis}, we sometimes use the following auxiliary continuous system:
\begin{align}
	\label{def EMCon}
	\begin{split}
		& Y_0=x,   \\
		& {\rm d} Y_{t}=b(Y_{t_{k-1}}) {\rm d} t+\sigma{\rm d} B_t,\ t\in(t_{k-1},t_{k}],\quad k \ge 1.
	\end{split}
\end{align} 
It is easy to see that $Y_{t_{k}}$ and $\theta_{k}$ have the same distribution for every $k$.

\subsection{Organization of this paper} 
In the next section, we first provide  assumptions and then present  main results, Theorems \ref{TheoremCLT} and \ref{TheoremFCLT}, and the byproducts of independent interest, Theorems \ref{TheoremW2distance} and \ref{TheoremCriterion}. We then provide a discussion of the results and point out that a special step-size sequence $\{1/k\}_{k \ge 1}$ may lead to the non-existence of a CLT, accompanied by simulation verification. The remainder of the paper is devoted to proving the CLT and the FCLT. To this end, we establish several supporting results; see Figure \ref{fig:map} below.
\begin{figure}[t]
	\centering
	\includegraphics[width=0.7\linewidth, height=0.5\linewidth]{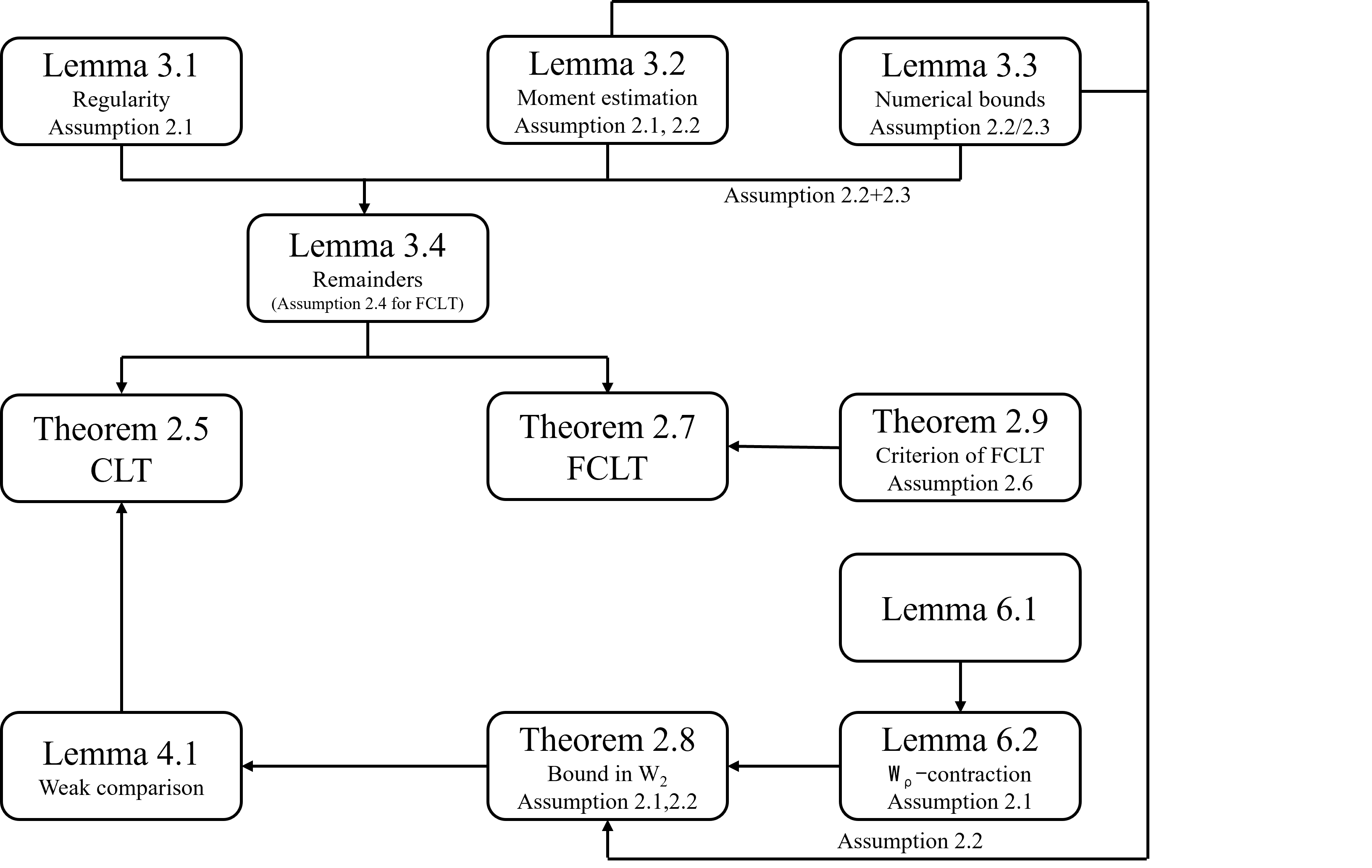}
	\caption{Proof structure}
	\label{fig:map}
\end{figure} 

In Section \ref{SectionPrelim}, we introduce decomposition \eqref{def Decomposition}, which is used to prove both the CLT and the FCLT. The remainder term in \eqref{def Decomposition} is justified to be negligible in Lemma \ref{LemRem}, proved by Lemmas \ref{LemmaRegularityStein}, \ref{LemMomentEstim} and \ref{LemNum}. Those lemmas are proved in Appendix \ref{SectionProofPrelim}.

With Lemma \ref{LemRem} in hand, it is sufficient to prove the CLT and the FCLT for the martingale part in \eqref{def Decomposition}. The CLT is proved in Section \ref{SectionProofCLT} by weak comparison in Lemma \ref{LemmaWeakComparison}. As a byproduct, we obtain a $\mathbf{W}_{2}$--distance error bound (Theorem \ref{TheoremW2distance}). Theorem \ref{TheoremW2distance} is proved in Section  \ref{SectionProofW2distance}, preceded by some other auxiliary lemmas proved in Appendix \ref{SectionProofLemmaThms}. The FCLT is proved in Section \ref{SectionProofIP} by verifying the   criterion  Theorem \ref{TheoremCriterion}. Since the criterion is   a slight modification of \ocite{McL74}*{Theorem 3.2},  we include a proof of Theorem \ref{TheoremCriterion} in Appendix \ref{AppendixProofIP} for completeness.
 
\section{Assumptions and Main results}\label{SectionMainResults}
\subsection{Assumptions}
We assume the following for coefficients of both \eqref{def SDE} and \eqref{def EMDis}.
\begin{assumption}
	\label{Asp Coeff}
	Assume that the function $b:\mathbb{R}^{d}\to\mathbb{R}^{d}$   is continuous and satisfies the following conditions:\\
	$\bullet$ $($Lipschitz condition$)$: There exists some constant $L>0$, such that \begin{align*}
		\vert  b(x)-  b(y)\vert\le L\vert x-y\vert,\quad \forall x,y\in\mathbb{R}^{d};
	\end{align*}
	$\bullet$ $($Partial dissipation$)$: There  exist constants $K_{1}>0$, $K_{2}\ge 0$, such that \begin{align*}
		\langle  b(x)-b(y),x-y\rangle\le -K_{1}\vert x-y\vert^{2}+K_{2},\quad \forall x,y\in\mathbb{R}^{d}.
	\end{align*}
	Assume that there exists $K_{3}>1$, such that the matrix $\sigma\in\mathbb{R}^{d\times d}$  is invertible and satisfies\begin{align*}
		K_{3}^{-1}I_{d}\preceq \sigma\sigma^{\intercal}\preceq K_{3}I_{d}.
	\end{align*}
\end{assumption}

Assumption \ref{Asp Coeff} implies that for some $C>0$, for all $x,y\in\mathbb{R}^{d}$,
\begin{align}
	\label{ine DispOneSide}
	\langle x, b(x)\rangle&\le-\frac{K_{1}}{2}\vert x\vert^{2}+C, \\
	\label{ine LinearGrowth}\vert b(x)\vert&\le C(1+\vert x\vert), 
\end{align}
and
\begin{align}
	\label{ine Boundsigmay}K_{3}^{-1}\vert y\vert^{2}\le \vert \sigma^{-1}y\vert^{2}\le K_{3}\vert y\vert^{2} .
\end{align}

It is easy to see that the infinitesimal generator of the SDE \eqref{def SDE} is \begin{align}\label{def generatorX}
	\mathcal{A}f(x)=\langle b(x),\nabla f(x)\rangle+\frac{1}{2}\langle \sigma\sigma^{\intercal},\nabla^{2}f(x)\rangle_{\text{HS}},
\end{align}
and \ocite{LTX22}*{Lemma 2.3} show that for $V(x)=1+\vert x\vert^{2}$, one has \begin{align}\label{ine Lyapnov}
	\mathcal{A}V(x)\le -\lambda V(x)+q
\end{align}
for some constants $\lambda,q>0$. 
\begin{assumption}\label{Asp Etak-0}
	The sequence of step sizes $\{\eta_{k}\}_{k\ge 1}\subset(0,1)$ satisfies:\\
	$\bullet$ $\{\eta_{k}\}_{k\ge 1}$ is non-increasing, i.e.  \begin{align*}
	 	\eta_{1} \ge \eta_{2} \ge \cdots \ge\eta_{k} \ge \eta_{k+1} \ge \cdots >0.
	\end{align*}
	$\bullet$ The partial sum $t_{n}$ defined in \eqref{def t_{n}} is divergent, i.e. \begin{align*}
		\sum_{k=1}^{\infty}\eta_{k}=\infty.
	\end{align*}
	$\bullet$ There exists some constant $c>0$, such that \begin{align}
		\label{con etak-diff ine}
		\eta_{k-1}-\eta_{k}\le c\eta_{k}^{2},\quad k\ge 2.
	\end{align} 
\end{assumption}
Assumptions \ref{Asp Coeff} and \ref{Asp Etak-0} give the moment estimations of \eqref{def SDE} and \eqref{def EMDis}; see Lemma \ref{LemMomentEstim} below. To obtain CLT, we need Assumption \ref{Asp Etak-1} below.

\begin{assumption}\label{Asp Etak-1}
	For the sequence  $\{\eta_{k}\}_{k\ge 1}\subset(0,1)$, let $T_{n}$ be defined as in \eqref{def T_n}. Assume: \\
	$\bullet$ For some $0<\epsilon<1$, 
	\begin{align}\label{con etak-strongR-M}
		\sum_{k=1}^{\infty}\eta_{k}^{2-\epsilon}<\infty.
	\end{align}
	$\bullet$  As $n\to\infty$, for some $L_{n}\to \infty$ sufficiently slowly, one has  
	\begin{align}\label{con etak-rate of etanTn}
		\frac{L_{n}}{\eta_{n}\sqrt{T_{n}}}\to 0.
	\end{align}   
\end{assumption}
To ensure FCLT, we impose additionally Assumption \ref{Asp Etak-2}. 
\begin{assumption}\label{Asp Etak-2}
 	For the sequence  $\{\eta_{k}\}_{k\ge 1}\subset(0,1)$, let $T_{n}$ be defined as in \eqref{def T_n}. There exists some $p>1$, such that as $n\to\infty$,  
 	\begin{align}\label{con etak-rate uniform}
 		\frac{(1+t_{n})^{1/p}}{\eta_{n}\sqrt{T_{n}}}\to 0.
 	\end{align}  
\end{assumption}
\begin{remark}~{}\\
	$\bullet$ The classical   \ocite{RM51} condition reads as: \begin{align*} 
		\sum_{k=1}^{\infty}\eta_{k}=\infty ,\quad  \sum_{k=1}^{\infty}\eta_{k}^{2}<\infty.
	\end{align*}
	However, due to the lower regularity of the solution to the Poisson equation \eqref{def Poisson equation}, we  require a  slightly stronger condition \eqref{con etak-strongR-M}.\\
	$\bullet$ The term $L_{n}$ in \eqref{con etak-rate of etanTn} is introduced for technical reasons, and  its choice is therefore not unique; see proof of Lemma \ref{LemNum} for details. The condition \eqref{con etak-rate of etanTn} implies that $L_{n}^{2}\lesssim \eta_{n}^{2}T_{n}$, but the proof only requires $L_{n}\lesssim \eta_{n}^{2}T_{n}\le n\eta_{n}\le t_{n}$; we keep the current form for consistency.\\
	$\bullet$ The condition \eqref{con etak-rate of etanTn} also implies that for any constant $C>0$, as $n\to\infty$, \begin{align}\label{con etak-rate weak}
		\frac{C}{\eta_{n}\sqrt{T_{n}}}\to 0.
	\end{align} 
	
\end{remark}
\subsection{The main results}
Consider $\Pi_{n}$ with a test function $h\in\text{\rm{Lip}}(\mathbb{R}^{d};\mathbb{R})$, i.e. \begin{align*}
	\Pi_{n}(h)=\frac{1}{n}\sum_{k=0}^{n-1}\delta_{\theta_{k}}(h)=\frac{1}{n}\sum_{k=0}^{n-1}h(\theta_{k}).
\end{align*}

\noindent Under assumptions above, we shall prove the following CLT.
\begin{theorem}\label{TheoremCLT}
	Suppose that Assumptions \ref{Asp Coeff}, \ref{Asp Etak-0} and \ref{Asp Etak-1} hold.   Assume that $\eta_1 \le \eta^{\ast}$ for some sufficiently small $\eta^{\ast}>0$. Let $h\in \text{\rm{Lip}}(\mathbb{R}^{d};\mathbb{R})$. Then we have\begin{align}
		\label{CLT}
		\frac{n}{\sqrt{T_{n}}}[\Pi_{n}(h)-\pi(h)] \Rightarrow \mathcal{N}(0,\pi(\vert \sigma^{\intercal}\nabla\varphi\vert^{2})),\quad\text{ as }n\to\infty,
	\end{align}
	where $\varphi$ solves Poisson equation \begin{align}
		\label{def Poisson equation}
		h-\pi(h)=\mathcal{A}\varphi,
	\end{align}
	and $\mathcal{A}$ is the generator of SDE \eqref{def SDE} given by \eqref{def generatorX}. 
\end{theorem}
\begin{remark}
	The order of $n/\sqrt{T_{n}}$ is determined by the choice of step sizes $\{\eta_{k}\}_{k\ge 1}$. In the classical CLT framework, the order of time scaling factor is $n^{1/2}$; see, for example, \ocite{Dur19}*{Theorem 3.4.1}. In our setting, $n/\sqrt{T_{n}}$ varies with $\{\eta_{k}\}_{k\ge1}$. For instance, when $\eta_{k}=k^{-\beta}$, $n/\sqrt{T_{n}}$ is of order $n^{(1-\beta)/2}$.
\end{remark} 
Furthermore, we establish an FCLT related to Theorem \ref{TheoremCLT}. The limiting process is a standard Brownian motion  time-changed by $a(t)$. For $T>0$, we partition $[0,T]$ into abutting intervals $[t_{\delta}^{m},t_{\delta}^{m+1}]$, each of length $\delta$, except possibly the final interval, whose length is at most $\delta$. We impose regularity conditions of $a(t)$ through its behavior at grid points $t_{\delta}^{m}$'s, as specified in the following assumption.
\begin{assumption}\label{Asp a}
	Let $T>0$, write $[0,T]=\bigcup_{m=0}^{\lfloor \delta^{-1}T\rfloor}[t_{\delta}^{m},t_{\delta}^{m+1}]$, where $t_{\delta}^{m}=m\delta$ for $0\le m\le \lfloor\delta^{-1}T\rfloor$ and $t_{\delta}^{\lfloor \delta^{-1}T\rfloor+1}=T$. Let $a(t):[0,T]\to\mathbb{R}_{+}$ be strictly increasing continuous function with $a(0)=0$. Set $\tau_{m}(\delta)\coloneqq a(t_{\delta}^{m+1})-a(t_{\delta}^{m})$ and $\tau(\delta)\coloneqq\sup_{0\le m\le \lfloor \delta^{-1}T\rfloor}\tau_{m}(\delta)$. Assume that \begin{align*}
		\lim_{\delta\to 0}\tau(\delta)\log(\delta^{-1})=0. 
	\end{align*}
\end{assumption}
 
\begin{theorem}\label{TheoremFCLT} 
	Let the assumptions of Theorem \ref{TheoremCLT} and Assumption \ref{Asp Etak-2} hold. Assume that the limit  \begin{align*}
		a(t)\coloneqq\lim_{n\to\infty}\frac{T_{\lfloor nt \rfloor}}{T_{n}}\in[0,1]
	\end{align*}
	exists and satisfies Assumption \ref{Asp a} with $T=1$. Then one has \begin{align*}
		\left\{\frac{\lfloor nt\rfloor }{\sqrt{T_{n}}}[\Pi_{\lfloor nt \rfloor}(h)-\pi(h)]\right\}_{t \in [0,1]}\Rightarrow \sqrt{\pi(\vert\sigma^{\intercal}\nabla\varphi\vert^{2})}\{B^{1}_{a(t)}\}_{t \in [0,1]}
	\end{align*}
	as $n\to\infty$ on $D[0,1]$, where $\{B^{1}_{t}\}_{t\ge 0}$ is one-dimensional standard Brownian motion and $\{B^{1}_{a(t)}\}_{t\in[0,1]}$ is time-changed Brownian motion.
\end{theorem}

\noindent \underline{\bf Two examples.} Consider $\eta_{k}=k^{-\alpha}$, $\alpha\in(\frac{1}{2-\epsilon},1)$ with $\epsilon$ as in Assumption \ref{Asp Etak-1}. It is easy to check that $a(t)=t^{1+\alpha}$, which satisfies Assumption \ref{Asp a}, so that Theorem \ref{TheoremFCLT} tells us 
\begin{align*}
	\left\{\frac{\lfloor nt\rfloor}{\sqrt{T_{n}}}[\Pi_{\lfloor nt \rfloor}(h)-\pi(h)]\right\}_{t \in [0,1]}\Rightarrow\sqrt{\pi(\vert\sigma^{\intercal}\nabla\varphi\vert^{2})}\left\{B_{t^{1+\alpha}}^{1}\right\}_{t \in [0,1]}.
\end{align*} 
Consider $\eta_{k}=\frac{\log(2+k)}{2+k}$, then $a(t)=t^{2}$, which satisfies Assumption \ref{Asp a}, by Theorem \ref{TheoremFCLT}, 
\begin{align*}
	\left\{\frac{\lfloor nt\rfloor}{\sqrt{T_{n}}}[\Pi_{\lfloor nt \rfloor}(h)-\pi(h)]\right\}_{t \in [0,1]}\Rightarrow\sqrt{\pi(\vert\sigma^{\intercal}\nabla\varphi\vert^{2})}\left\{B_{t^{2}}^{1}\right\}_{t \in [0,1]}.
\end{align*} 
\subsection{By-products of independent interest}
When proving Theorem \ref{TheoremCLT}, we need the crucial Lemma \ref{LemmaWeakComparison} below, whose proof heavily depends on a bound of $\mathbf{W}_{2}(X_{t_{n}},\theta_{n})$. Such bound may be of independent interest.   
\begin{theorem}\label{TheoremW2distance} 
	Suppose that Assumption \ref{Asp Coeff} and Assumption \ref{Asp Etak-0} hold.  Assume that $\eta_{1} \le \eta^{\ast}$ for some sufficiently small $\eta^{\ast}>0$. Let $\{X_{t}\}_{t\ge0}$ and $\{\theta_{k}\}_{k\ge0}$ be given as in \eqref{def SDE} and \eqref{def EMCon} respectively. For all $n\ge 1$, any $0\le i<n$, with $X_{t_{i}}=\theta_{i}=z\in\mathbb{R}^{d}$ be given, we have
	\begin{align*}
		\mathbf{W}_{2}(X_{t_{n}},\theta_{n})\le C(1+\vert z\vert^{1/2})\eta_{n}^{1/4}.
	\end{align*}  
\end{theorem}
\begin{remark}
	\ocite{MMS20} proved the $\mathbf{W}_{2}$--distance convergence rate for ULA to the target distribution $\pi$ with a fixed step size is of order $\eta^{1/4}$. Their approach employs a novel coupling technique, combining synchronous coupling and the coupling in \ocite{EM19}. \ocite{MMS20}*{Section 5.1} briefly discussed the possible extension to the case of decreasing step sizes.
\end{remark}
When proving Theorem \ref{TheoremFCLT}, the classical criterion  \ocite{McL74}*{Theorem 3.2} is not applied directly. Following the argument in its proof, we make a tiny generalization of this criterion below, which immediately yields our FCLT. We put its proof in Appendix \ref{AppendixProofIP}.
\begin{theorem}\label{TheoremCriterion}
	Let $\{k_{n}(t)\}_{n\ge 1}$ be a sequence of integer-valued, non-decreasing, and right-continuous functions defined on $[0,T]$ with $k_{n}(0)=0$ for all $n\ge 1$. Suppose $X_{n,i}$ is a martingale difference array satisfying that for each $t\in[0,T]$, as $n\to\infty$:\\
	$(1)$ $\max_{i\le k_{n}(t)}\vert X_{n,i}\vert\to0$ in $L^{2}$;\\
	$(2)$ $\sum_{i=1}^{k_{n}(t)}X_{n,i}^{2}\to a(t)$ in probability.\\
	Let \begin{align*}
		S_{n}(t)\coloneqq\sum_{i=1}^{k_{n}(t)}X_{n,i},\ t\in[0,T].
	\end{align*}
	Then, under Assumption \ref{Asp a}, as $n\to\infty$,  \begin{align*}
		\{S_{n}(t)\}_{t\in[0,T]}\Rightarrow \{B_{a(t)}^{1}\}_{t\in[0,T]}
	\end{align*}
	on $D[0,T]$, where $B^{1}$ denotes one-dimensional standard Brownian motion.
\end{theorem} 
\subsection{Discussions on a special case}
Let us discuss a special step-size sequence $\eta_{k}=1/k$ with $k \ge 1$, we have  \begin{align*}
	\frac{1}{\eta_{n}\sqrt{T_{n}}}\asymp 1,
\end{align*}  so that even the weak form \eqref{con etak-rate weak} of \eqref{con etak-rate of etanTn} in Assumption \ref{Asp Etak-1} is not satisfied. We conduct simulations to investigate whether the CLT fails in this case.

We consider the one dimensional case, setting $\sigma=1$ and choosing
\begin{align*}
	b(x)=-x+\sin x.
\end{align*}

\emph{Experiment 1}:  Let the Lipschitz test function be 
\begin{align*}
	h(x)=\frac{1}{1+x^{2}}.
\end{align*}
Run the EM scheme $N=3000$ times, each with $n=10^{5}$ iterations whose initializations are uniformly chosen from $\{-8,2,12\}$. We consider 
\begin{align*}
	\eta_{k}=\frac{1}{k^{3/4}},\quad  \text{ and }\quad \eta_{k}=\frac{1}{k}.
\end{align*}

It is easy to check that the former choice satisfies Assumptions \ref{Asp Etak-0}, \ref{Asp Etak-1}, and \ref{Asp Etak-2}, whereas the latter does not. We simulate the empirical distributions of \eqref{def prelim sum}.

\begin{figure}[t]\caption{Result of \emph{Experiment 1}}\label{Fig1}	\centering
	\subfigure[$\eta_{k}=\frac{1}{k^{3/4}}$]{\begin{minipage}[t]{0.45\linewidth} 
			\includegraphics[width=0.8\linewidth]{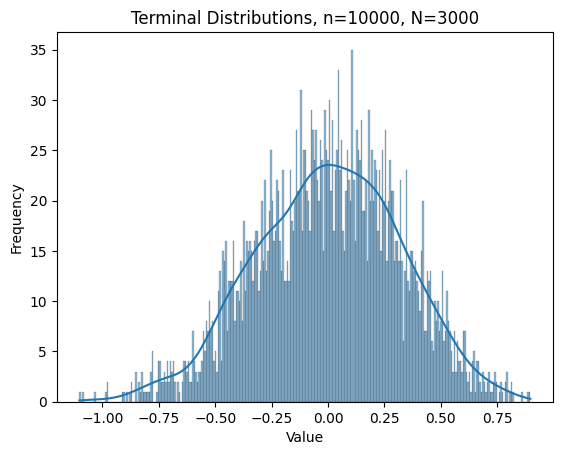} 
	\end{minipage}}
	\subfigure[$\eta_{k}=\frac{1}{k}$]{\begin{minipage}[t]{0.45\linewidth} 
			\includegraphics[width=0.8\linewidth]{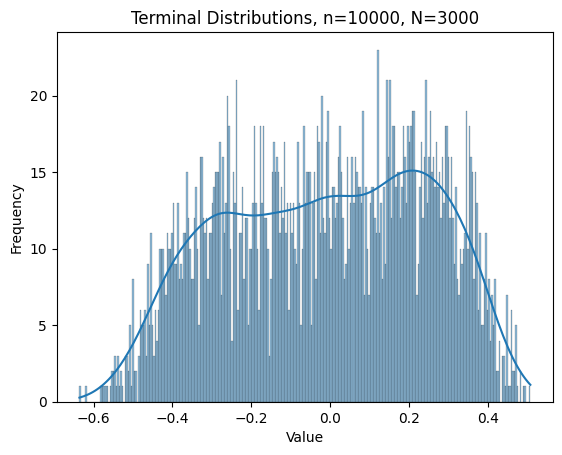} 
	\end{minipage}}
\end{figure}

The results are shown in Figure \ref{Fig1}, which suggests that the CLT may not hold when Assumption \ref{Asp Etak-1} is violated. 

One may doubt whether such failure is due to insufficient number of iterations. To test this concern, we fix the step size sequence as $\eta_{k}=1/k$, and increase the number of iterations.

\emph{Experiment 2}: Let the Lipschitz test functions be \begin{align*}
	h(x)=\frac{1}{1+x^{2}},\quad  h(x)=\sin x,\quad h(x)=x.
\end{align*} 
Run the EM scheme $N=3000$ times, each with $n=10^{9}$ iterations whose initializations are uniformly chosen from $\{-8,2,12\}$. The step sizes are $\eta_{k}=1/k$, $k\ge 1$. 

\begin{figure}[t]
	\centering\caption{Result of \emph{Experiment 2}}\label{Fig2}
	\subfigure[$h(x)=\frac{1}{1+x^{2}}$]{\begin{minipage}[t]{0.3\linewidth} 
			\includegraphics[width=0.8\linewidth]{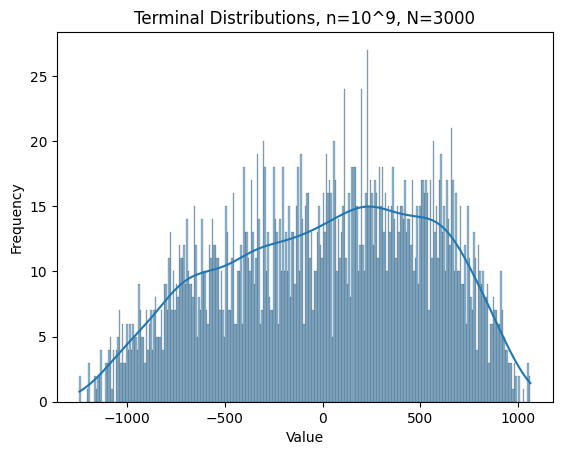} 
	\end{minipage}}
	\subfigure[$h(x)=\sin x$]{\begin{minipage}[t]{0.3\linewidth} 
			\includegraphics[width=0.8\linewidth]{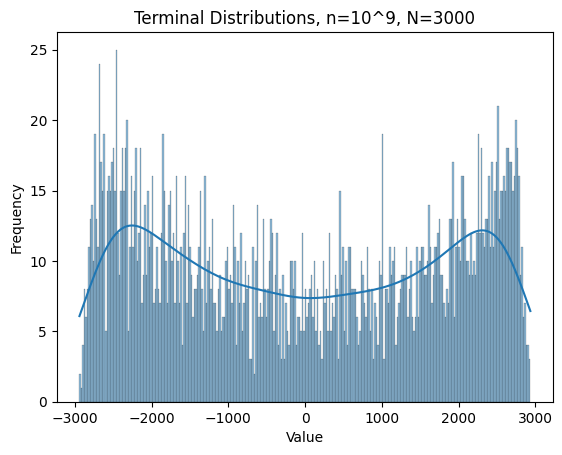} 
	\end{minipage}}
	\subfigure[$h(x)=x$]{\begin{minipage}[t]{0.3\linewidth} 
			\includegraphics[width=0.8\linewidth]{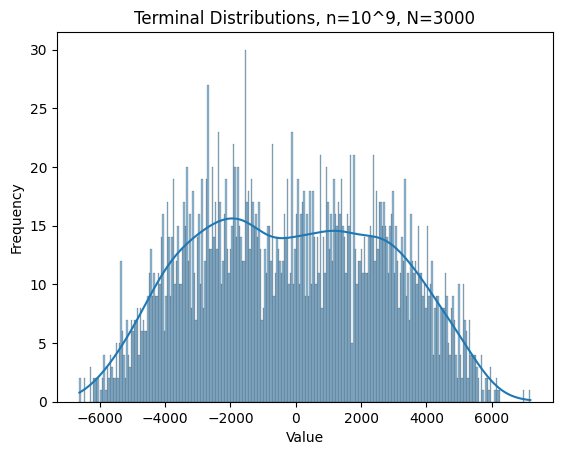} 
	\end{minipage}}	
\end{figure}

The results are shown in Figure \ref{Fig2}, which shows significant variation across different test functions, suggesting that the CLT does not hold when $\eta_{k}=1/k$. We leave it as a problem.
\begin{problem}For the step-size sequence  $\eta_{k}=1/k$, $k\ge 1$, let $T_{n}$ be defined   \eqref{def T_n}. Let $\{\theta_{k}\}_{k\ge0}$ be as in \eqref{def EMDis}. Let $h\in \text{\rm Lip}(\mathbb{R}^{d};\mathbb{R})$, and let $\Pi_{n}$ be the empirical measure defined in \eqref{def Pi_n}.  Determine the limiting distribution of \eqref{def prelim sum} as $n\to\infty$. 
\end{problem}
 
\section{Decomposition, Poisson equation and auxiliary lemmas}\label{SectionPrelim}
Let us give the strategy for proving our main results, which includes a decomposition, a Poisson equation, and auxiliary lemmas for the proofs. 

\subsection{The strategy for proving CLT} Let $\varphi: \mathbb{R}^d \rightarrow \mathbb R$ be a function with a certain regularity such that the following Taylor expansion holds:   
\begin{align*} 
	\varphi(\theta_{k+1})-\varphi(\theta_{k})=\langle\nabla\varphi(\theta_{k}),\Delta\theta_{k+1}\rangle+\frac{1}{2}\langle\nabla^{2}\varphi(\theta_{k}),(\Delta\theta_{k+1})(\Delta\theta_{k+1})^{\intercal}\rangle_{\text{HS}}+\mathcal{R}_{k+1}, 
\end{align*} 
where\begin{align}\label{defDeltaThetak}
	\Delta\theta_{k+1}\coloneqq\theta_{k+1}-\theta_{k}=\eta_{k+1}b(\theta_{k})+\sqrt{\eta_{k+1}}\sigma\xi_{k+1},
\end{align}
\begin{align}\label{def R_k+1}
	\mathcal{R}_{k+1}\coloneqq\int_{0}^{1}\int_{0}^{1}r\langle\nabla^{2}\varphi(\theta_{k}+r  s  \Delta\theta_{k+1}) -\nabla^{2}\varphi(\theta_{k}),(\Delta\theta_{k+1})(\Delta\theta_{k+1})^{\intercal}\rangle_{\text{HS}}{\rm d}r{\rm d}s.
\end{align}
Recalling the definition of $\mathcal{A}$ from \eqref{def generatorX}, we have 
\begin{align} \label{equation difference phi-phi-Aphi}
	\begin{split}
		&\varphi(\theta_{k+1})-\varphi(\theta_{k})-\eta_{k+1}\mathcal{A} \varphi(\theta_k) \\ 
		=&\sqrt{\eta_{k+1}}\langle\nabla\varphi(\theta_{k}),\sigma \xi_{k+1}\rangle+\frac{1}{2}\langle\nabla^{2}\varphi(\theta_{k}),(\Delta\theta_{k+1})(\Delta\theta_{k+1})^{\intercal}-\eta_{k+1}
		\sigma \sigma^{\intercal}\rangle_{\text{HS}}+\mathcal{R}_{k+1}.
	\end{split}
\end{align} 
In order to make a connection between the above expansion and our CLT, we introduce the Poisson equation \eqref{def Poisson equation}, i.e., 
\begin{align*}
	\mathcal{A}\varphi=h-\pi(h).
\end{align*}
Let $\varphi$ be the solution to the Poisson equation. For an integer $m\le n$, we obtain:
\begin{align*}
	\sum_{k=0}^{m-1}[h(\theta_{k})-\pi(h)]
	&=\sum_{k=0}^{m-1}\frac{1}{\eta_{k+1}}\mathcal{A}\varphi(\theta_{k})\eta_{k+1}\\
	&=\sum_{k=0}^{m-1}\frac{1}{\eta_{k+1}}\left[\mathcal{A}\varphi(\theta_{k})\eta_{k+1}-[\varphi(\theta_{k+1})-\varphi(\theta_{k})]\right]+\sum_{k=0}^{m-1}\frac{\varphi(\theta_{k+1})-\varphi(\theta_{k})}{\eta_{k+1}}.
\end{align*}
By \eqref{equation difference phi-phi-Aphi}, we get the following decomposition:  
\begin{align}
	\label{def Decomposition}
	\frac{1}{\sqrt{T_{n}}}\sum_{k=0}^{m-1}[h(\theta_{k})-\pi(h)]=\frac{1}{\sqrt{T_{n}}}M_{m}+R_{n}(m),
\end{align}
where
\begin{align}\label{def M_m Z_m}
	M_{m}=\sum_{k=1}^{m}Z_{k},\quad Z_{k}\coloneqq -\frac{1}{\sqrt{\eta_{k}}}\langle \nabla\varphi(\theta_{k-1}),\sigma\xi_{k}\rangle,
\end{align}
and \begin{align}\label{def R_n(m)}
	R_{n}(m)=R_{n}^{(0)}(m)+R_{n}^{(1)}(m)-R_{n}^{(2)}(m)-R_{n}^{(3)}(m).
\end{align}
We give the  forms of remainders before proving their negligibility; see Section \ref{SubsecProofLemRem} below.

To prove CLT, we apply  \eqref{def Decomposition} with $m=n$, use the CLT criterion for the martingale $M_{m}/\sqrt{T_{n}}$, and prove that $R_{n}(n)$ is negligible in probability. To prove FCLT, we apply \eqref{def Decomposition} with $m=\lfloor nt\rfloor$, use our FCLT criterion Theorem \ref{TheoremCriterion} for $M_{\lfloor nt\rfloor}/\sqrt{T_{n}}$, and prove that $\sup_{t\in[0,1]}R_{n}(\lfloor nt\rfloor)$ is negligible in probability.
 
To advance this strategy, we require the following auxiliary lemmas, which encompass the regularity of the Poisson equation, the moment estimates associated with SDE and the bounds with respect to step sizes.
\subsection{Poisson equation \eqref{def Poisson equation} and regularity of its solution}
Let $\varphi$ be the solution to Poisson equation \eqref{def Poisson equation}, namely: 
\begin{align*}
	\mathcal{A}\varphi=h-\pi(h).
\end{align*}
When the test function $h\in C_{b}^{2}(\mathbb{R}^{d};\mathbb{R})$, it follows that the solution $\varphi\in C_{b}^{4}$ with $\Vert\nabla^{k}\varphi\Vert_{\infty}\le C$, $k=0,1,2,3,4$, as shown in \ocite{LTX22}*{Lemma 3.1}. In this paper, we  consider  $h\in \text{\rm{Lip}}(\mathbb{R}^{d};\mathbb{R})$  instead, the corresponding regularity of $\varphi$ is given by the following lemma.
\begin{lemma}[\ocite{DFL26}*{Lemma 3.3}]\label{LemmaRegularityStein}
	Suppose Assumption \ref{Asp Coeff} holds. Let $h\in\text{\rm{Lip}}(\mathbb{R}^{d};\mathbb{R})$, and let $\varphi$ be solution to Poisson equation \eqref{def Poisson equation}. Then there exists a constant $C>0$, such that \begin{align}\label{bound |grad_k phi(x)|}
		\vert \nabla^{k}\varphi(x)\vert\le C(1+\vert x\vert^{k+2}),\ k=0,1,2.
	\end{align}
	Further, it follows that \begin{align}\label{bound DiffHesPhi}
		\sup_{y:\vert y-x\vert\le 1}\frac{\vert\nabla^{2}\varphi(x)-\nabla^{2}\varphi(y)\vert}{\vert x-y\vert}\le C(1+\vert x\vert^{5}).
	\end{align}
\end{lemma}

\subsection{Remainder terms in decomposition \eqref{def Decomposition} are  negligible}
We give the two auxiliary lemmas about properties of  $\{\theta_{k}\}_{k\ge0}$ and $\{\eta_{k}\}_{k\ge1}$; see Appendix \ref{SectionProofPrelim} for their proofs. 
\begin{lemma} \label{LemMomentEstim}
	Let $\{X_{t}\}_{t\ge0}$ be the solution to SDE \eqref{def SDE}, and $\{\theta_{k}\}_{k\ge 1}$ be the associated EM scheme defined by \eqref{def EMDis}. Then, under  Assumption \ref{Asp Coeff}, one has:\\
	$\bullet$ For any $p\ge 1$, there exists a constant $C_{1}>0$ such that \begin{align}\label{bound momentEstimX}
		\sup_{t\ge s}\mathbb{E}[\vert X_{t}\vert^{p}|X_s=z]\le C_{1}(1+\vert z \vert^{p})
	\end{align}
	for all $s \ge 0$ and $z \in \mathbb{R}^d$.\\ 
	Under Assumptions \ref{Asp Coeff} and \ref{Asp Etak-0}, assume further that there exists $\eta^{\ast}>0$ sufficiently small, such that $\eta_{1}\le \eta^{\ast}$. Then:\\
	$\bullet$ For any $p\ge 1$, there exists a   constant $C_{2}>0$ such that for all $i \ge 0$ and $z \in \mathbb{R}^d$. 
	\begin{align} \label{bound E|b(theta_k)|^p}
		\sup_{k \ge i} \mathbb{E}[\vert  b(\theta_{k})\vert^{p}|\theta_i=z]\le C_{2}(1+\vert z \vert^{p}), \quad \sup_{k \ge i} \mathbb{E}[\vert \theta_{k}\vert^{p}|\theta_i=z]\le C_{2}(1+\vert z \vert^{p}).
	\end{align}
	$\bullet$ For any $p\ge 1$, there exists a  constant $C_{3}>0$, depending on $p$, such that \begin{align}\label{bound theta_k^2p uniformly}
		\mathbb{E}\left[\max_{0\le k\le n}(1+\vert \theta_{k}\vert^{2})^{p}\right]\le C_{3}(1+t_{n})(1+\vert x\vert^{2p}).
	\end{align}
\end{lemma}
\begin{lemma}\label{LemNum}
	Under Assumption \ref{Asp Etak-0}, as $n\to\infty$, the following statements hold:\\
	$\bullet$
	\begin{align}
		\label{LimNum1}
		T_{n}^{-2}\sum_{1\le i<j\le n}\eta_{i}^{-3/4}\eta_{j}^{-1}\to 0.
	\end{align}
	$\bullet$  For any $n\ge 2$, any constant $C_{1}>0$, and any $q\in[0,1]$, there exists a constant $C_{2}>0$ independent of $n$, such that: \begin{align}\label{bound SumEtaExp Num2}
		\sum_{k=2}^{n}e^{-C_{1}(t_{n}-t_{k})}\eta_{k}^{1+q}\le C_{2}\eta_{n}^{q}.
	\end{align}
		Under  Assumptions \ref{Asp Etak-0} and \ref{Asp Etak-1}, as $n\to\infty$, the following statements hold:\\ 
	$\bullet$ For any $p\ge1$, \begin{align}
		\label{LimNum3}
		T_{n}^{-1/2}\sum_{k=1}^{n}\eta_{k}^{p/2}\to 0.
	\end{align} 
	$\bullet$  For any constant $C_{3}>0$, \begin{align}\label{LimNum4}
		T_{n}^{-2}\sum_{k=1}^{n}\frac{1}{\eta_{k}^{2}}\to 0,\quad T_{n}^{-2} \sum_{1\le i<j\le n}\eta_{i}^{-1}\eta_{j}^{-1}e^{-C_{3}(t_{j}-t_{i})}\to 0.
	\end{align}
\end{lemma}
Combining lemmas above, we obtain the following result:
\begin{lemma}\label{LemRem}
		Under  Assumptions \ref{Asp Coeff}, \ref{Asp Etak-0} and \ref{Asp Etak-1},  \begin{align*}
		\vert R_{n}^{(i)}(n)\vert\to 0\quad\text{ in probability},\quad i=0,1,2,3.
	\end{align*} 
	Under   Assumptions \ref{Asp Coeff}, \ref{Asp Etak-0}, \ref{Asp Etak-1}, and \ref{Asp Etak-2}, one has
	\begin{align*}
		\sup_{t\in[0,1]}\vert R_{n}^{(i)}(\lfloor nt\rfloor)\vert\to 0\quad\text{ in probability}, \quad i=0,1,2,3.
	\end{align*}
\end{lemma}
\begin{remark}\label{RemarkRn0}~{}\\
	$\bullet$ The error term $R_{n}^{(0)}(n)$ will vanish automatically in the homogeneous Markov chain setting by stationary initialization. However, in our inhomogeneous case, we need to use the Abel transform to show its negligibility, which requires \eqref{LimNum4} from condition \eqref{con etak-rate of etanTn}.\\
	$\bullet$ Due to the worse regularity of $\varphi$, we cannot further expand the term $R_{n}^{(3)}(m)$. To prove its negligibility, we need the condition $\eqref{con etak-strongR-M}$ in Assumption \ref{Asp Etak-1}, which gives \eqref{LimNum3}.\\
	$\bullet$ The FCLT requires uniform control of $\sup_{t\in[0,1]}R_{n}^{(1)}(\lfloor nt\rfloor)$, and the  condition $\eqref{con etak-rate of etanTn}$ is insufficient in such case. We additionally impose \eqref{con etak-rate uniform} in Assumption \ref{Asp Etak-2}.
\end{remark}

\section{Proof of Theorem \ref{TheoremCLT} (CLT)}\label{SectionProofCLT}
As stated before, we prove CLT by applying the martingale CLT \ocite{McL74}*{Theorem 2.3} on $M_{n}/\sqrt{T_{n}}$. To this end, we need the following error bound between $\vert\sigma^{\intercal}\nabla\varphi(\theta_{k})\vert^{2}$ and $\pi(\vert\sigma^{\intercal}\nabla\varphi\vert^{2})$; see Appendix \ref{SectionProofLemmaThms} below for its proof.
\begin{lemma}
	\label{LemmaWeakComparison}Under the assumptions of Theorem \ref{TheoremCLT}, for $i<j-1$,  \begin{align}
		\label{WeakComparisonGiven i}
		\left\vert\mathbb{E}_{i}[\vert\sigma^{\intercal}\nabla\varphi(\theta_{j-1})\vert^{2}-\pi(\vert\sigma^{\intercal}\nabla\varphi\vert^{2})]\right\vert \le C_{1}V_{6}(\theta_{i})e^{-C_{0}(t_{j}-t_{i})}+C_{2}(1+\vert \theta_{i}\vert^{15/2})\eta_{j-1}^{1/4},
	\end{align}
	holds for some constants $C_{0},C_{1},C_{2}>0$, where $\mathbb{E}_{i}[\cdot]\coloneqq\mathbb{E}[\cdot\vert\theta_{i}]$, $V_{6}(x)=1+\vert x\vert^{6}$, $t_{i}$ and $t_{j}$ are defined in \eqref{def t_{n}}, and $\varphi$ solves Poisson equation \eqref{def Poisson equation}.
\end{lemma} 
 
\begin{proof}[Proof of Theorem \ref{TheoremCLT}]
	Recall the decomposition \eqref{def Decomposition}. Combining Lemma \ref{LemRem} and the well-known Slutsky's Theorem  \ocite{Slu25}, it is sufficient to show that as $n\to\infty$,
	\begin{align} \label{LimW MartCLT}
		\frac{1}{\sqrt{T_{n}}}M_{n} \Rightarrow N(0,\pi(\vert\sigma^{\intercal}\nabla \varphi\vert^{2})).
	\end{align}
	Recall notations \eqref{def M_m Z_m}. According to \ocite{McL74}*{Theorem 2.3}, \eqref{LimW MartCLT} holds if we can verify that, as $n\to\infty$, the following conditions are true:\begin{align}\tag{C1}
		\label{McLCond1}\mathbb{E}\max_{1\le i\le n} \frac{\vert Z_{i}\vert}{\sqrt{T_{n}}}\to 0
	\end{align}
	and \begin{align}\tag{C2}
		\label{McLCond2}\sum_{k=1}^{n} \frac{Z_{k}^{2}}{T_{n}}\to\pi(\vert\sigma^{\intercal}\nabla\varphi\vert^{2})\quad\text{ in probability}.
	\end{align} 
	
	\underline{Verification of \eqref{McLCond1}:} We turn to verify
	\begin{align}\label{C1prime}
		\tag{C1'}\max_{1\le i\le n} \frac{\vert Z_{i}\vert}{\sqrt{T_{n}}}\to 0\quad\text{ in }L^{2},
	\end{align}
	which implies \eqref{McLCond1} by H\"{o}lder's inequality. It will also be applied in proving FCLT. To verify \eqref{C1prime}, for each $1\le i\le n$, splitting by $A_{i}\coloneqq\{\vert Z_{i}\vert^{2}\le \eta_{n}^{-2}\}$ yields
	\begin{align}\label{C1Split}
		\left(\max_{1\le i\le n}\frac{\vert Z_{i}\vert}{\sqrt{T_{n}}}\right)^{2}=\frac{1}{T_{n}}\max_{1\le i\le n}\vert Z_{i}\vert^{2}=\frac{1}{T_{n}}\max_{1\le i\le n}\vert Z_{i} 1_{A_{i}}\vert^{2}+\frac{1}{T_{n}}\max_{1\le i\le n}\vert Z_{i} 1_{A_{i}^{c}}\vert^{2}.
	\end{align} 
	The first term in \eqref{C1Split} is bounded by  \eqref{con etak-rate of etanTn}:  \begin{align}\label{McLCond1Step1}
		\mathbb{E}\frac{1}{T_{n}}\max_{1\le i\le n}\vert Z_{i} 1_{A_{i}}\vert^{2}\le \frac{1}{\left(\eta_{n}\sqrt{T_{n}}\right)^{2}}\to 0.
	\end{align}
	For bounding second term in \eqref{C1Split}, by   \eqref{bound |grad_k phi(x)|}, \eqref{bound E|b(theta_k)|^p}, \eqref{ine Boundsigmay} and the independence of $\xi_{i}$ and $\theta_{i-1}$, we have for each $1\le i\le n$,
	\begin{align*}
		\mathbb{E}\vert Z_{i}\vert^{2}\le \frac{1}{\eta_{i}}\mathbb{E}\vert \nabla\varphi(\theta_{i-1})\vert^{2}\mathbb{E}\vert\sigma\xi_{i}\vert^{2} \le \frac{C}{\eta_{i}}(1+\mathbb{E}\vert\theta_{i-1}\vert^{6})\mathbb{E}\vert\xi_{i}\vert^{2}\le \frac{C}{\eta_{i}}.
	\end{align*}
	Similarly, $\mathbb{E}\vert Z_{i}\vert^{4}\le C\eta_{i}^{-2}$. Hence,
	\begin{align*}
		\mathbb{E}[\vert Z_{i}\vert^{2} 1_{A_{i}^{c}}]\le[\mathbb{E}\vert Z_{i}\vert^{4}]^{1/2}[\mathbb{P}(A_{i}^{c})]^{1/2}\le \frac{C}{\eta_{i}}\sqrt{\mathbb{P}\left(\vert Z_{i}\vert^{2}>\eta_{n}^{-2}\right)} \le\frac{C\eta_{n}}{\eta_{i}^{3/2}},
	\end{align*} 
	where the last inequality follows from Chebyshev's inequality. By monotonicity, one has  $\eta_{i}^{-1/2}\le \eta_{n}^{-1/2}$ for all $1\le i\le n$, canceling a $T_{n}$ yields\begin{align}\label{McLCond1Step2}
		\mathbb{E}\frac{1}{T_{n}}\max_{1\le i\le n}\vert Z_{i}1_{A_{i}^{c}}\vert^{2}\le\frac{1}{T_{n}}\sum_{i=1}^{n}\mathbb{E}[\vert Z_{i}\vert^{2} 1_{A_{i}^{c}}]\le \frac{C\eta_{n}}{T_{n}}\sum_{i=1}^{n}\frac{1}{\eta_{i}^{3/2}}\le C\eta_{n}^{1/2}\to 0
	\end{align} 
	as $n\to\infty$. Combine \eqref{McLCond1Step1} and \eqref{McLCond1Step2} to get \eqref{C1prime}. 
	
	\underline{Verification of \eqref{McLCond2}:} It suffices to show the following $L^2$ convergence: as $n\to\infty$, 
	\begin{align*}
		\mathbb{E}\left[\sum_{k=1}^{n}\frac{Z_{k}^{2}}{T_{n}}-\pi(\vert\sigma^{\intercal}\nabla\varphi\vert^{2})\right]^{2}\to0. 
	\end{align*}
	By expanding the square directly and taking absolutes value, one has
	\begin{align}	\label{bound C2L2} 
			\mathbb{E}\left[\sum_{k=1}^{n}\frac{Z_{k}^{2}}{T_{n}}-\pi(\vert\sigma^{\intercal}\nabla\varphi\vert^{2})\right]^{2}\le\frac{1}{T_{n}^{2}}\sum_{k=1}^{n}\frac{1}{\eta_{k}^{2}}\mathbb{E} J_{k}+\frac{1}{T_{n}^{2}}\sum_{1\le i<j\le n}\frac{1}{\eta_{i}}\frac{1}{\eta_{j}}\vert \mathbb{E} J_{ij}\vert, 
	\end{align} 
	where 
	\begin{align*}
		J_{k}&\coloneqq[\langle\nabla\varphi(\theta_{k-1}),\sigma\xi_{k}\rangle^{2}-\pi(\vert\sigma^{\intercal}\nabla\varphi\vert^{2})]^{2},\\
		J_{ij}&\coloneqq[\langle\nabla\varphi(\theta_{i-1}),\sigma\xi_{i}\rangle^{2}-\pi(\vert\sigma^{\intercal}\nabla\varphi\vert^{2})][\langle\nabla\varphi(\theta_{j-1}),\sigma\xi_{j}\rangle^{2}-\pi(\vert\sigma^{\intercal}\nabla\varphi\vert^{2})],
	\end{align*}
	and we estimate them below.

	By the Cauchy--Schwartz inequality, \eqref{bound |grad_k phi(x)|}, \eqref{ine Boundsigmay}, and  moment estimation \eqref{bound E|b(theta_k)|^p}, we have\begin{align*}
		\mathbb{E} J_{k}\le \mathbb{E}\vert\nabla\varphi(\theta_{k-1})\vert^{4}\mathbb{E}\vert\sigma\xi_{k}\vert^{4} +2\pi(\vert\sigma^{\intercal}\nabla\varphi\vert^{2})^{4}\le C(1+\vert x\vert^{12}),
	\end{align*} 
	which gives\begin{align*}
		\vert \mathbb{E} J_{ij}\vert\le[\mathbb{E} J_{i}]^{1/2}[\mathbb{E} J_{j}]^{1/2}\le C(1+\vert x\vert^{12}).	
	\end{align*} 
	Therefore, combining monotonicity, we have
	\begin{align}\label{bound C2L2Quad}
			 \sum_{k=1}^{n}\frac{1}{\eta_{k}^{2}}\mathbb{E} J_{k}+\sum_{\substack{1\le i<j\le n,\\ i=j-1}}\frac{1}{\eta_{i}}\frac{1}{\eta_{j}}\vert \mathbb{E} J_{ij}\vert    \le  C (1+\vert x\vert^{12})\sum_{k=1}^{n}\frac{1}{\eta_{k}^{2}}. 
	\end{align}
	Let us estimate the remaining terms $J_{ij}$ with $j-1>i$. Conditioning on $\theta_{i}$ gives \begin{align*}
		&\left\vert \mathbb{E} J_{ij}\right\vert\\
		\le & \mathbb{E}\left[ \left\vert \langle\nabla\varphi(\theta_{i-1}),\sigma\xi_{i}\rangle^{2}-\pi(\vert\sigma^{\intercal}\nabla\varphi\vert^{2})  \right\vert\left\vert\mathbb{E}_{i}[\vert\sigma^{\intercal}\nabla\varphi(\theta_{j-1})\vert^{2}-\pi(\vert\sigma^{\intercal}\nabla\varphi\vert^{2})]\right\vert \right] \\
		\le & \left[ \mathbb{E}\vert \nabla\varphi(\theta_{i-1})\vert^{2} \vert \sigma\xi_{i}\vert^{2}+\vert \pi(\vert\sigma^{\intercal}\nabla\varphi\vert^{2})\vert^{2}\right]^{1/2} \left[\mathbb{E}\left\vert\mathbb{E}_{i}[\vert\sigma^{\intercal}\nabla\varphi(\theta_{j-1})\vert^{2}-\pi(\vert\sigma^{\intercal}\nabla\varphi\vert^{2})]\right\vert^{2}\right]^{1/2},
	\end{align*} 
	where $\mathbb{E}_{i}[\cdot]\coloneqq\mathbb{E}[\cdot\vert \theta_{i}]$. 
	By \eqref{bound |grad_k phi(x)|} and \eqref{bound E|b(theta_k)|^p}, we have 
	\begin{align*}
		\mathbb{E}\vert\nabla\varphi(\theta_{i-1})\vert^{2} \vert \sigma\xi_{i}\vert^{2} +\vert \pi(\vert\sigma^{\intercal}\nabla\varphi\vert^{2})\vert^{2} \le C(1+\vert x\vert^{6}). 
	\end{align*}
	First applying Lemma \ref{LemmaWeakComparison} to $\left\vert\mathbb{E}_{i}[\vert\sigma^{\intercal}\nabla\varphi(\theta_{j-1})\vert^{2}-\pi(\vert\sigma^{\intercal}\nabla\varphi\vert^{2})]\right\vert$ and then using \eqref{bound |grad_k phi(x)|} and \eqref{bound E|b(theta_k)|^p}, we have 
	\begin{align*}
		\mathbb{E}\left\vert\mathbb{E}_{i}[\vert\sigma^{\intercal}\nabla\varphi(\theta_{j-1})\vert^{2}-\pi(\vert\sigma^{\intercal}\nabla\varphi\vert^{2})]\right\vert^{2} \le C e^{-2C_{0}(t_{j-1}-t_i)}(1+\vert x\vert ^{12})+C(1+\vert x\vert^{15}) \eta^{1/2}_{j-1}.
	\end{align*}
	Therefore, for all $j-1>i$, using $e^{-C_{0}(t_{j-1}-t_{i})} =e^{-C_{0}(t_{j-1}-t_{j})} e^{-C_{0}(t_{j}-t_{i})} \le e^{C_{0}\eta_{1}} e^{-C_{0}(t_{j}-t_{i})}$ with noting $\eta_{i} \ge \eta_{j-1}$,  one has
	\begin{align}\label{bound C2L2Cross}
		\vert \mathbb{E} J_{ij}\vert \le C(1+\vert x\vert ^{11}) \left[e^{-C_{0}(t_{j-1}-t_{i})}+\eta_{j-1}^{1/4}\right]   \le  C(1+\vert x\vert ^{11}) \left[e^{-C_{0}(t_{j}-t_{i})}+\eta_{i}^{1/4}\right].
	\end{align}

	Combining \eqref{bound C2L2}, \eqref{bound C2L2Quad} and \eqref{bound C2L2Cross}, and using Lemma \ref{LemNum}, we have as $n\to\infty$,
	\begin{align*}
		&\mathbb{E}\left[\sum_{k=1}^{n}\frac{Z_{k}^{2}}{T_{n}}-\pi(\vert\sigma^{\intercal}\nabla\varphi\vert^{2})\right]^{2}\\
		\le&\frac{C(1+\vert x\vert ^{12})}{T_{n}^{2}}\left[\sum_{1\le i<j\le n}\frac{1}{\eta_{i}}\frac{1}{\eta_{j}}e^{-C_{0}(t_{j}-t_{i})}+\sum_{1\le i<j\le n }\frac{1}{\eta_{i}^{3/4}}\frac{1}{\eta_{j}}+\sum_{k=1}^{n}\frac{1}{\eta_{k}^{2}}\right]\to0.
	\end{align*}
	The condition \eqref{McLCond2} is verified.
\end{proof}
 
\section{Proof of Theorem \ref{TheoremFCLT} (FCLT)}\label{SectionProofIP}
With Lemma \ref{LemRem} in hand, it is sufficient to prove the FCLT for $M_{\lfloor nt\rfloor}/\sqrt{T_{n}}$ by verifying conditions of criterion Theorem \ref{TheoremCriterion}. 
\begin{proof}[Proof of Theorem \ref{TheoremFCLT}]
	Recall the decomposition \eqref{def Decomposition}, since Lemma \ref{LemRem} shows that as $n\to\infty$ its remainder part \begin{align*}
		\sup_{t\in[0,1]}\vert R_{n}(\lfloor nt\rfloor)\vert\to 0\quad \text{ in probability}, 
	\end{align*}it is sufficient to prove the FCLT for the martingale part of \eqref{def Decomposition}:\begin{align*}
		\frac{1}{\sqrt{T_{n}}}M_{\lfloor nt\rfloor}=\sum_{k=1}^{\lfloor nt\rfloor}\frac{Z_{k}}{\sqrt{T_{n}}}.
	\end{align*}
	To this end, we verify conditions in Theorem \ref{TheoremCriterion}. Recall that in proving Theorem \ref{TheoremCLT}, the condition $(1)$ of Theorem \ref{TheoremCriterion} has been checked as in \eqref{C1prime}, and by \eqref{McLCond2} one has 
	\begin{align*}
		\frac{1}{\pi(\vert\sigma^{\intercal}\nabla\varphi\vert^{2})}\sum_{k=1}^{n}\frac{Z_{k}^{2}}{T_{n}}\to 1\quad \text{ in }L^{2},\quad\text{ as } n\to\infty,
	\end{align*}
	which implies for each $t\in[0,1]$ that \begin{align*}
		\frac{1}{\pi(\vert\sigma^{\intercal}\nabla\varphi\vert^{2})}\sum_{k=1}^{\lfloor nt \rfloor}\frac{Z_{k}^{2}}{T_{n}}=\frac{1}{\pi(\vert\sigma^{\intercal}\nabla\varphi\vert^{2})}\frac{T_{\lfloor nt \rfloor}}{T_{n}}\sum_{k=1}^{\lfloor nt \rfloor}\frac{Z_{k}^{2}}{T_{\lfloor nt \rfloor}}\to a(t)\quad \text{ in }L^{2}
	\end{align*}
	as $n\to\infty$. By Theorem \ref{TheoremCriterion}, the result follows.
\end{proof}
 
\section{Proof of Theorem \ref{TheoremW2distance} ($\mathbf{W}_{2}$--convergence)}\label{SectionProofW2distance}
\subsection{The framework and some notations for coupling}\label{SubsecCoupling}
We shall use the reflection coupling framework in \ocite{EGZ19} to prove Theorem \ref{TheoremW2distance}. Since we need to bound the $\mathbf{W}_{2}$--distance ${\rm W_{2}}(X_{t_{k}},\theta_{k})$, the result therein cannot be directly applied. Instead, we use the reflection coupling method in \ocites{Ebe16,EGZ19} for each time interval $[t_{i}, t_{i+1}]$ for $0 \le i \le k-1$.

Let us briefly introduce the notation in \ocite{EGZ19}. We shall use their  $\rho_{1}$--distance: \begin{align}\label{def rho1}
	\rho_{1}(x,y)\coloneqq 1_{\{x\ne y\}}[f(\vert x-y\vert)+\varepsilon V(x)+\varepsilon V(y)],
\end{align}
where $f$ is a non-decreasing concave continuous function satisfying $f(0)=0$, $\varepsilon>0$ is a positive constant to be chosen later, and $V$ is a Lyapunov function that will be chosen in our setting as 
\begin{align*}
	V(x)=1+\vert x\vert^{2}. 
\end{align*}
It is known that (see for example  \ocite{Vil09}*{Theorem 6.15}) the $\mathbf{W}_{2}$--distance between probability measures $\mu$ and $\nu$ is controlled by a weighted total variation norm 
\begin{align*}
	\mathbf{W}_{2}(\mu,\nu)\le \sqrt{2}\left(\int\vert z\vert^{2}\vert\mu-\nu\vert({\rm d}z)\right)^{1/2}.
\end{align*}
Combine this with the facts $\vert z\vert^{2}< 1+\vert z\vert^{2}=V(z)$ for $z\in\mathbb{R}^{d}$ and \begin{align*}
	\int_{\mathbb{R}^{d}}V(z)\vert\mu-\nu\vert (dz)=\inf_{\gamma\in \textbf{C}(\mu,\nu)}\int[V(x)+V(y)]1_{\{x\ne y\}}\gamma({\rm d}x,{\rm d}y)
\end{align*}
as shown in \ocite{Hai11}*{Lemma 2.1}, it follows that \begin{align}\label{ControlW2ByWrho1}
	\mathbf{W}_{2}(\mu,\nu)\le (2/\varepsilon)^{-1/2}[\mathbf{W}_{\rho_{1}}(\mu,\nu)]^{1/2}.
\end{align}

For completeness, let us recall the construction of $f$ in \ocite{Ebe16}, which will also be used in our setting. By Assumption \ref{Asp Coeff}, there exists a function $\kappa:[0,\infty)\to\mathbb{R}$ satisfying 
\begin{align}
	\label{DissipativeDeduced}
	\langle  b(x)-b(y),x-y\rangle\le \kappa(\vert x-y\vert)\vert x-y\vert^{2},
\end{align}
\begin{align}
	\lim_{r\to\infty}\kappa(r)=-K_{1}<0,\quad \lim_{\delta\to 0}\sup_{r\in[0,\delta]}r\kappa(r)= 0,\quad \int_{0}^{\infty}r[\kappa(r)\lor 0]{\rm d}r<\infty.
\end{align}
For example, one may take $\kappa(r)=\min\{-K_{1}+\frac{K_{2}}{r^{2}},L\}$. For such $\kappa$, define $R_{0}\coloneqq\inf\{s\ge 0:\kappa(r)\le 0,\ \forall r\ge s\}$ and $R_{1}\coloneqq\inf\{s\ge R_{0}:s(s-R_{0})\kappa(r)\le - 8 ,\ \forall r\ge s\}$. Further, define function $\varphi,\Phi:\mathbb{R}_{+}\to\mathbb{R}_{+}$ by \begin{align*}
	\varphi(r)\coloneqq\exp\left(- \frac{K_{3}}{2}\int_{0}^{r}s[\kappa(s)\lor 0]{\rm d}s\right), \Phi(r)=\int_{0}^{r}\varphi(s){\rm d}s,
\end{align*}
and the function $g$ is defined as
\begin{align*}
	g(r)\coloneqq1-\frac{c_{1}K_{3}}{2}\int_{0}^{r\land R_{1}}\Phi(s)\varphi^{-1}(s){\rm d}s-c_{2}K_{3}\int_{0}^{r\land R_{1}}\varphi(s)^{-1}{\rm d}s,
\end{align*}
where constants $c_{1}$ and $c_{2}$ are given by
\begin{align}\label{c1c2}
	c_{1}\coloneqq\left[2K_{3}\int_{0}^{R_{1}}\Phi(r)\varphi(r)^{-1}{\rm d}r\right]^{-1},\ c_{2}\coloneqq\left[4K_{3}\int_{0}^{R_{1}}\varphi(r)^{-1}{\rm d}r\right]^{-1}.
\end{align} 
Using $\varphi$ and $g$, we can define a non-decreasing concave function $f:\mathbb{R}_{+}\to \mathbb{R}_{+}$ by setting \begin{align}\label{def  f}
	f(r)\coloneqq\int_{0}^{r}\varphi(s)g(s){\rm d}s. 
\end{align} 
It can be verified that $f$ is a twice differentiable function, such that for any $r\ge0$, \begin{align}\label{boundedfprime}
	0\le f^{\prime}(r)\le 1,\quad \frac{\varphi(R_{0})}{2}r\le f(r)\le r
\end{align}
and for $r\in[0,R_{1}]$, 
\begin{align}
	\label{fDoublePrime}
	\frac{1}{K_{3}}f^{\prime\prime}(r)+ \frac{1}{2} r\kappa(r)f^{\prime}(r)\le -\frac{c_{1}}{2}f(r)-c_{2}.
\end{align}
Further, one can check  for $r\ge 0$ that 
\begin{align}
	\label{fDoubleprimeGeneral}
	\frac{1}{K_{3}}f^{\prime\prime}(r)+\frac{1}{2}r\kappa(r)f^{\prime}(r)\le-\frac{1}{2}c_{1}^{\prime}f(r)
\end{align} 
for some $c_{1}^{\prime}>0$; cf. \ocite{Ebe16}*{Section 4}.

For fixed $\delta>0$, let $\phi_{1}^{\delta},\phi_{2}^{\delta}:\mathbb{R}^{d}\to[0,1]$ be two continuous and Lipschitz functions satisfying \begin{align}\label{Phi1+Phi2=1}
	[\phi_{1}^{\delta}(x)]^{2}+[\phi_{2}^{\delta}(x)]^{2}=1,\ \forall x\in\mathbb{R}^{d},\quad\text{ and }\quad\phi_{1}^{\delta}(x)=\left\{\begin{array}{ll}
		1,&\vert x\vert\ge\delta,\\
		0,&\vert x\vert\le\frac{\delta}{2}.
	\end{array}\right.
\end{align}
By \eqref{boundedfprime} and the definition of $\phi_{2}^{\delta}$, it is easy to see that for all $x\in\mathbb{R}^{d}$,
\begin{align}\label{ThreeUpperBounds}
	[\phi_{2}^{\delta}(x)]^{2}\le 1_{\{\vert x\vert<\delta\}},\quad \ [\phi_{2}^{\delta}(x)]^{2}f(\vert x\vert)\le \delta,\quad \ [\phi_{2}^{\delta}(x)]^{2}\vert x\vert\kappa(\vert x\vert)f^{\prime}(\vert x\vert)\le \sup_{r\in[0,\delta]}r\kappa(r).
\end{align}
\subsection{The coupling argument and auxiliary lemmas}
We will work with the continuous systems. Given $X_{t_i}=z$ and $Y_{t_i}=z$, let us consider the two stochastic processes $\{X_{t}\}_{t_{i} \le t \le t_{n}}$ and $\{Y_{t}\}_{t_{i} \le t \le t_{n}}$ which are correspondingly governed by 
\begin{align}
	\label{SDEcoupling}{\rm d} X_{t}=b(X_{t}) {\rm d}t+\sigma{\rm d} B_{t},\quad X_{t_{i}}=z;
\end{align} 
and 
\begin{align}
	\label{EMcoupling} \begin{split}
		& Y_{t_{i}}=z,  \\
		& {\rm d} Y_{t}=b(Y_{t_{k-1}}) {\rm d} t+\sigma{\rm d} B_{t},\ t\in(t_{k-1},t_{k}],\quad  i+1 \le k \le n. 
	\end{split}
\end{align}
To bound $\mathbf{W}_{2}$--distance,  according to \eqref{ControlW2ByWrho1}, we instead seek a $W_{\rho_{1}}$-distance upper bound. By definition of $W_{\rho_{1}}$-distance, the problem reduces to bound $\mathbb{E}\rho_{1}(X_{t_{n}},Y_{t_{n}})$. To this end, we shall follow exactly the reflection coupling method developed in \ocite{Ebe16}.

Consider the following reflection coupling for $t\in[t_{n-1},t_{n}]$, $n-1\ge i$:
\begin{align}
\label{Coupling}
\begin{split}
	{\rm d}\widetilde{X}_{t}&=b(\widetilde{X}_{t}){\rm d}t+\sigma\left[\phi_{1}^{\delta}(Z_{t}) {\rm d}B_{t}^{(1)}+\phi_{2}^{\delta}(Z_{t}){\rm d}B_{t}^{(2)}\right],\\
	{\rm d}\widetilde{Y}_{t}&=b(\widetilde{Y}_{t_{n-1}}){\rm d}t+\phi_{1}^{\delta}(Z_{t})\sigma\left[I_{d}-\frac{2(\sigma^{-1}z_{t})(\sigma^{-1}z_{t})^{\intercal}}{\vert\sigma^{-1} z_{t}\vert^{2}}\right]{\rm d}B_{t}^{(1)}+\phi_{2}^{\delta}(Z_{t})\sigma{\rm d}B_{t}^{(2)},
\end{split}
\end{align}
where $\widetilde{X}_{t_i}=\widetilde{Y}_{t_i}=z$,
\begin{align}
	\label{def Z,z,n}Z_{t}\coloneqq\widetilde{X}_{t}-\widetilde{Y}_{t},\quad z_{t}\coloneqq n(Z_{t}),\quad n(x)\coloneqq\left\{\begin{array}{ll}
		\frac{x}{\vert x\vert},&x\in\mathbb{R}^{d}\setminus\{0\},\\
		0,&x=0.
	\end{array}\right.
\end{align}
and $\{B_{t}^{(1)}\}_{t\ge t_i}$, $\{B_{t}^{(2)}\}_{t\ge t_i}$ are two mutually independent standard Brownian motions. Under  Assumption \ref{Asp Coeff}, it is easy to see that $\widetilde{X}_{t}\stackrel{d}{=}X_{t}$ and $\widetilde{Y}_{t}\stackrel{d}{=}Y_{t}$.
\begin{lemma}\label{LemmaDifferenceCoupling}
	Let $\{\widetilde{X}_{t}\}_{t_{n-1}\le t\le t_{n}}$, $\{\widetilde{Y}_{t}\}_{t_{n-1}\le t\le t_{n}}$ be defined as in \eqref{Coupling}. Let $Z_{t}$ be defined as in \eqref{def Z,z,n}. Let $i\ge 0$, for all $n\ge i+1$, $t\in[t_{n-1},t_{n}]$, one has\begin{align}
		{\rm d}\vert Z_{t}\vert= \langle z_{t}, b(\widetilde{X}_{t})-b(\widetilde{Y}_{t_{n-1}})\rangle {\rm d}t+2\phi_{1}^{\delta}(Z_{t})\frac{ (\sigma^{-1}z_{t})^{\intercal}}{\vert \sigma^{-1}z_{t}\vert^{2}}{\rm d}B_{t}^{(1)}.
	\end{align}
\end{lemma}
We see from the definition \eqref{def rho1} of $\rho_{1}$ that 
\begin{align}
	\label{dErho(X,Y)}
	\frac{\rm d}{\rm dt}\mathbb{E} \rho_{1}(X_{t},Y_{t})\le \frac{\rm d}{\rm dt}\mathbb{E} f(\vert Z_{t}\vert)+ \frac{\rm d}{\rm dt}\mathbb{E} [\varepsilon V(X_{t})+\varepsilon V(Y_{t})], 
\end{align} 
where the function $f$ is chosen as in \eqref{def  f}. Starting from \eqref{dErho(X,Y)}, we can derive the following result, which yields an upper bound of $\mathbb{E} \rho_{1}(X_{t_{n}},Y_{t_{n}})$ as desired.
\begin{lemma}\label{LemmaContraction}Let $\{X_{t}\}_{t_{n-1}\le t\le t_{n}}$ and $\{Y_{t}\}_{t_{n-1}\le t\le t_{n}}$ be defined as in \eqref{SDEcoupling} and \eqref{EMcoupling} respectively.  There exist constants $C_{1},C_{2}>0$, such that for every sufficiently small $\delta>0$ and $t\in[t_{n-1},t_{n}]$,
	\begin{align}\label{ddtErho1}
		\frac{{\rm d}}{{\rm d}t}\mathbb{E} \rho_{1}(X_{t},Y_{t}) 
		\le-C_{1}\mathbb{E} \rho_1(X_{t},Y_{t})+c_{1}\delta+\sup_{r\in[0,\delta]}r\kappa(r)+2c_{2}\mathbb{P}\left(\vert Z_{t}\vert<\delta\right)+C_{2}(1+\vert z\vert)\eta_{n}^{1/2},
	\end{align}
	where constants $c_{1},c_{2}>0$ are defined as in \eqref{c1c2}.
\end{lemma}
We then use lemmas above to prove Theorem \ref{TheoremW2distance}. The proofs of these lemmas can be found in Appendix \ref{SectionProofLemmaThms}.
\subsection{Proof of Theorem \ref{TheoremW2distance}}
\begin{proof}[Proof of Theorem \ref{TheoremW2distance}]
	Given $X_{t_i}=z$ and $Y_{t_i}=z$. Consider the stochastic processes $\{X_{t}\}_{t_{i} \le t \le t_{n}}$ and $\{Y_{t}\}_{t_{i} \le t \le t_{n}}$ defined in \eqref{SDEcoupling} and \eqref{EMcoupling}. 	We shall follow the reflection coupling method developed in \ocite{EGZ19} to show 
	\begin{align} \label{bound Rho by sum}
		\mathbb{E}\rho_{1}(X_{t_{n}},Y_{t_{n}})\le C(1+\vert z\vert)\sum_{k=i+1}^{n}e^{-c(t_{n}-t_{k})}\eta_{k}^{3/2},
	\end{align} 
	which, together with the definition of ${\rm W}_{\rho_1}$, implies 
	\begin{align*} 
		\mathbf{W}_{\rho_{1}}(\text{\rm{Law}}(X_{t_{n}}),\text{\rm{Law}}(Y_{t_{n}}))\le C(1+\vert z\vert)\sum_{k=i+1}^{n}e^{-c(t_{n}-t_{k})}\eta_{k}^{3/2}. 
	\end{align*}
	By \eqref{ControlW2ByWrho1} and Lemma \ref{LemNum}, we get the bound 
	\begin{align*}
		\mathbf{W}_{2}(\text{\rm{Law}}(X_{t_{n}}),\text{\rm{Law}}(Y_{t_{n}}))\le C(1+\vert z\vert^{1/2}) \eta_{n}^{1/4}.
	\end{align*} 
	
	It remains to prove \eqref{bound Rho by sum}.  We shall follow exactly the method developed in \ocite{EGZ19}. Consider the processes $\{\widetilde{X}_{t}\}_{ t_{n-1}\le t\le t_{n}}$, $\{\widetilde{Y}_{t}\}_{ t_{n-1}\le t\le t_{n}}$ defined via reflection coupling \eqref{Coupling} for $t\in[t_{n-1},t_{n}]$,	where $\widetilde{X}_{t_i}=\widetilde{Y}_{t_i}=z$, $Z_{t}\coloneqq\widetilde{X}_{t}-\widetilde{Y}_{t}$. From Lemma \ref{LemmaContraction}, multiplying $e^{-C_{1}t}$ on both sides of \eqref{ddtErho1}, integrating over $[t_{n-1},t_{n}]$ and letting $\delta\to0$, we arrive at \begin{align}\label{ExpectedRho1Step}
		\mathbb{E}\rho_{1}(X_{t_{n}},Y_{t_{n}})\le e^{-C_{1}\eta_{n}}\mathbb{E}\rho_{1}(X_{t_{n-1}},Y_{t_{n-1}})+C(1+\vert z\vert)\eta_{n}^{3/2}.
	\end{align}
	By iteration,  \eqref{bound Rho by sum} follows immediately.
	
\end{proof}
 
\appendix
\section{Proof of lemmas in Section \ref{SectionPrelim}}\label{SectionProofPrelim}
\subsection{Proof of Lemma \ref{LemMomentEstim}}
\begin{proof}[Proof of Lemma \ref{LemMomentEstim}]
	By Jensen's inequality we only need to consider the case $p\ge 2$.  Note by \eqref{ine DispOneSide} and Young's inequality that for all $x\in\mathbb{R}^{d}$, there exist constants $C_{1},C_{2}>0$ such that 
	\begin{align}\label{MoemntEstimCommon}
		\begin{split}
			&p\vert x\vert^{p-2}\langle x,b(x)\rangle+\frac{1}{2}p(p-1)\vert x\vert^{p-2}\vert\sigma\vert^{2}\\
			\le&-K_{2}p\vert x\vert^{p}+K_{1}p\vert x\vert^{p-2}+\frac{1}{2}p(p-1)\vert x\vert^{p-2}\vert \sigma\vert^{2}\\
			\le& C_{1}-C_{2}\vert x\vert^{p}.
		\end{split}
	\end{align} 
	
	\underline{Proof of \eqref{bound momentEstimX}} :  By It\^{o}'s formula and \eqref{MoemntEstimCommon}, there exists a martingale $\{M_{t}\}_{t\ge s}$ independent of $X_{s}$, such that for all $t\ge s$,
	\begin{align*}
		{\rm d}\vert X_{t}\vert^{p}-{\rm d}M_{t}\le&\left(p\vert X_{t}\vert^{p-2}\langle X_{t},b(X_{t})\rangle +\frac{1}{2}p(p-1)\vert X_{t}\vert^{p-2}\vert\sigma\vert^{2}\right){\rm d}t\\
		\le&\left(C_{1}-C_{2}\vert X_{t}\vert^{p}\right){\rm d}t,
	\end{align*}
	which implies for all $t\ge s$ that \begin{align*}
		\frac{{\rm d}}{{\rm d}t}\mathbb{E}[\vert X_{t}\vert^{p}\vert X_{s}=z]\le (C_{1}-C_{2}\mathbb{E}[\vert X_{t}\vert^{p}\vert X_{s}=z]),
	\end{align*}
	integrating over $[s,t]$ then yields\begin{align*}
		\mathbb{E}[\vert X_{t}\vert^{p}\vert X_{s}=z]\le e^{-C_{2}(t-s)}\mathbb{E}[\vert X_{s}\vert^{p}\vert X_{s}=z]+C_{1}\int_{s}^{t}e^{-C_{2}(t-r)}{\rm d}r\le \frac{C_{1}}{C_{2}}+\vert z\vert^{p}.
	\end{align*} 
	
	\underline{Proof of  \eqref{bound E|b(theta_k)|^p}}: According to \eqref{ine LinearGrowth}, we only need to show the second assertion by considering separately the range of $k\ge i$. To this end, we use continuous system \eqref{def EMCon} with recalling $Y_{t_{k}}\stackrel{d}{=}\theta_{k}$ for all $k$. First we note for all $m\ge 0$, $s\in[t_{m},t_{m+1}]$, given $Y_{t_{m}}$, one has 
	\begin{align*}
		Y_{s}-Y_{t_{m}}\stackrel{d}{=}(s-t_{m})b(Y_{t_{m}})+\sqrt{s-t_{m}}\sigma\xi_{m+1},
	\end{align*}
	where $\xi_{m+1}\sim\mathcal{N}(0,I_{d})$ is independent of $Y_{t_{m}}$, it implies that\begin{align}\label{WithinOneStep}
		\mathbb{E}\left[\vert Y_{s}-Y_{t_{m}}\vert^{p}\vert Y_{t_{m}}\right]\le C_{1}\eta_{m+1}^{p}(1+\vert Y_{t_{m}}\vert^{p})+C_{2}\eta_{m+1}^{p/2}.
	\end{align} 
	Note that \eqref{bound E|b(theta_k)|^p} holds trivially for $k=i$; by using independence of $\theta_{i}$ and $\xi_{i+1}$ with taking $m=i$, $s=t_{i+1}$ in \eqref{WithinOneStep}, the case $k=i+1$ follows. It remains to  consider the case when $k>i+1$. By It\^{o}'s formula and \eqref{MoemntEstimCommon}, for all $t\in[t_{k},t_{k+1}]$, there exist a martingale $\{M_{t}\}_{t\ge s}$ independent of $\theta_{i}$, such that 
	\begin{align*}
		{\rm d}\vert Y_{t}\vert^{p}-{\rm d}M_{t}&\le\left(p\vert Y_{t}\vert^{p-2}\left[\langle Y_{t},b(Y_{t})\rangle+\langle Y_{t},b(Y_{t_{k}})-b(Y_{t})\rangle\right]+\frac{1}{2}p(p-1)\vert Y_{t}\vert^{p-2}\vert\sigma\vert^{2}\right){\rm d}t\\
		&\le \left(C_{1}(1+\vert Y_{t}-Y_{t_{k}}\vert^{p})-C_{2}\vert Y_{t}\vert^{p}\right){\rm d}t,
	\end{align*}
	where the last step follows from \eqref{MoemntEstimCommon} and the Young's inequality. This implies that\begin{align}
		\nonumber&\mathbb{E}[\vert \theta_{k+1}\vert^{p}\vert \theta_{i}=z]\\
		\nonumber\le&e^{-C_{1}\eta_{k+1}}\mathbb{E}[\vert Y_{t_{k}}\vert^{p}\vert\theta_{i}=z]+C_{1}\eta_{k+1}+C_{1}\int_{t_{k}}^{t_{k+1}}\mathbb{E}[\vert Y_{s}-Y_{t_{k}}\vert^{p}\vert\theta_{i}=z]{\rm d}s\\
		\label{MomentEstimDisOneStep}\le&(C_{1}\eta_{k+1}^{1+p}+e^{-C_{2}\eta_{k+1}})\mathbb{E}[\vert Y_{t_{k}}\vert^{p}\vert\theta_{i}=z]+C_{3}\eta_{k+1},
	\end{align} 
	where the following estimation from independence of $\theta_{i}$ and $\xi_{k+1}$ for all $k\ge i$ and \eqref{WithinOneStep} has been applied in the last step:
	\begin{align*}
		\int_{t_{k}}^{t_{k+1}}\mathbb{E}[\vert Y_{s}-Y_{t_{k}}\vert^{p}\vert\theta_{i}=z]{\rm d}s\le C_{1}\eta_{k+1}^{1+p}\mathbb{E}[\vert Y_{t_{k}}\vert^{p}\vert\theta_{i}=z]+C_{2}\eta_{k+1}^{1+p/2}
	\end{align*}
	By recursively using \eqref{MomentEstimDisOneStep}, one has for $k>i+1$, \begin{align}\label{LastButOneStepMomentEstim}
		\mathbb{E}[\vert \theta_{k}\vert^{p}\vert\theta_{i}=z]\le C_{0}\mathbb{E}[\vert \theta_{i}\vert^{p}\vert\theta_{i}=z]+C\left[\eta_{k}+\sum_{m=i+1}^{k-1}\eta_{m}\prod_{j=m+1}^{k}[C_{2}\eta_{m}^{1+p}+e^{-C_{1}\eta_{m}}]\right].
	\end{align}
	By assumptions, one may choose $\eta^{\ast}$ sufficiently small,  such that for all $k\ge 1$, \begin{align*}
		0\le C_{2}\eta_{k}^{1+p}+e^{-C_{1}\eta_{k}}\le 1-\frac{1}{2}C_{1}\eta_{k}.
	\end{align*}
	Using $1-x\le e^{-x}$,  the sum in \eqref{LastButOneStepMomentEstim} can then be bounded by using \eqref{bound SumEtaExp Num2} with $q=0$. \eqref{bound E|b(theta_k)|^p} then follows from \eqref{LastButOneStepMomentEstim}.

	\underline{Proof of \eqref{bound theta_k^2p uniformly}}: We first work on each time interval and then union the bounds. Denote $I^{(m)}\coloneqq\left\{k\ge 0:m\le t_{k}< m+1\right\}$, $0\le m\le \lfloor t_{n}\rfloor$, then $I^{(m)}$'s are non-empty and disjoint. For $0\le k\le n$,  one can then find exact one $0\le m\le \lfloor t_{n}\rfloor$ such that $k\in I^{(m)}$, thus  $\{0,1,\cdots,n\}\subset\bigcup_{m=0}^{\lfloor t_{n}\rfloor}I^{(m)}$. We  discuss on each unit interval then union the bounds.
	
	Fix $m$, denote $\ell_{m}\coloneqq \min I^{(m)}$, $u_{m}\coloneqq \max I^{(m)}$, then $t_{u_{m}}-t_{\ell_{m}}\le 1$. For $k\in I^{(m)}$, write \begin{align*}
		\theta_{k}=\theta_{\ell_{m}}+\sum_{\ell=\ell_{m}+1}^{k}\eta_{\ell}b(\theta_{\ell-1})+\sum_{\ell=\ell_{m}+1}^{k}\sqrt{\eta_{\ell}}\sigma\xi_{\ell}.
	\end{align*}
	Denote \begin{align*}
		M_{j}^{(m)}\coloneqq \sum_{\ell=\ell_{m}+1}^{j}\sqrt{\eta_{\ell}}\sigma\xi_{\ell},\ j\in I^{(m)},\quad  A_{m}\coloneqq \vert\theta_{\ell_{m}}\vert+c+\max_{j\in I^{(m)}}\vert M_{j}^{(m)}\vert,
	\end{align*}
	where $c$ comes from \eqref{ine LinearGrowth},
	then one has \begin{align*}
		\vert\theta_{k}\vert\le A_{m}+c\sum_{\ell=\ell_{m}+1}^{k}\eta_{\ell}\vert\theta_{\ell-1}\vert\le A_{m}+c\sum_{\ell=\ell_{m}+1}^{k}\eta_{\ell}\max_{\ell_{m}\le r\le \ell-1}\vert\theta_{r}\vert,
	\end{align*}
	yielding that \begin{align}\label{bound theta_k pre-ite}
			\max_{\ell_{m}\le r\le k}\vert\theta_{r}\vert \le A_{m}+c\sum_{\ell=\ell_{m}+1}^{k}\eta_{\ell}\max_{\ell_{m}\le r\le \ell-1}\vert\theta_{r}\vert.
	\end{align}
	We claim that \begin{align}\label{bound max theta_k by iteration}
		\max_{\ell_{m}\le r\le k}\vert\theta_{r}\vert\le A_{m}\prod_{\ell=\ell_{m}+1}^{k}(1+c\eta_{\ell}).
	\end{align}
	We first use this claim to prove \eqref{bound theta_k^2p uniformly},  and then verify the claim. Using $1+x\le e^{x}$ and $\sum_{\ell=\ell_{m}+1}^{u_{m}}\eta_{\ell}\le t_{u_{m}}-t_{\ell_{m}}\le 1$, \eqref{bound max theta_k by iteration} then implies \begin{align*}
		\max_{\ell_{m}\le r\le u_{m}}\vert\theta_{r}\vert\le e^{c\sum_{\ell=\ell_{m}+1}^{u_{m}}\eta_{\ell}} A_{m}\le e^{c}A_{m},
	\end{align*}
	which gives \begin{align*}
		\mathbb{E}\max_{k\in I^{(m)}}\vert\theta_{k}\vert^{2p}=\mathbb{E}\max_{\ell_{m}\le r\le u_{m}}\vert\theta_{r}\vert^{2p}\le e^{2pc}\mathbb{E}A_{m}^{2p}\le C_{p}\left(1+\mathbb{E}\vert \theta_{\ell_{m}}\vert^{2p}+\mathbb{E}\max_{j\in I^{(m)}}\vert M_{j}^{(m)}\vert^{2p}\right).
	\end{align*}
	We turn to bound the last term. By \eqref{bound E|b(theta_k)|^p}, \begin{align*}
		\sup_{m\ge 0}\mathbb{E}\vert \theta_{\ell_{m}}\vert^{2p}\le C(1+\vert x\vert^{2p})
	\end{align*}
	for some $C>0$ independent of $m$; meanwhile, we note that the martingale $M_{j}^{(m)}\sim \mathcal{N}(0,\sigma\sigma^{\intercal}\sum_{\ell=\ell_{m}+1}^{j}\eta_{\ell})$ has uniformly in $m$ bounded covariance, by Doob's $L^{2p}$-inequality, \begin{align*}
		\mathbb{E}\max_{j\in I^{(m)}}\vert M_{j}^{(m)}\vert^{2p}\le C_{p}\mathbb{E} \vert M_{u_{m}}^{(m)}\vert^{2p}\le C
	\end{align*}
	for some $C$ independent of $m$. Therefore, one has \begin{align*}
		\mathbb{E}\max_{k\in I^{(m)}}\vert\theta_{k}\vert^{2p}\le C_{p}(1+\vert x\vert^{2p}). 
	\end{align*}
	Union the bounds, we get \begin{align*}
		\mathbb{E}\left[\max_{0\le k\le n}(1+\vert \theta_{k}\vert^{2})^{p}\right]\le\sum_{m=0}^{\lfloor t_{n}\rfloor}\mathbb{E}\max_{k\in I^{(m)}}(1+\vert \theta_{k}\vert^{2})^{p}\le C_{p}(1+t_{n})(1+\vert x\vert^{2p}).
	\end{align*}  
	
	It remains to verify \eqref{bound max theta_k by iteration}. The case $k=\ell_{m}$ is trivial; for $k=\ell_{m}+1$, \eqref{bound theta_k pre-ite} directly gives\begin{align*}
		\max_{\ell_{m}\le r\le \ell_{m}+1}\vert\theta_{r}\vert\le A_{m}+c\eta_{\ell_{m}+1}\vert \theta_{\ell_{m}}\vert\le A_{m}(1+c\eta_{\ell_{m}+1});
	\end{align*} 
	now assume \eqref{bound max theta_k by iteration} holds until $k-1$ for $k> \ell_{m}+1$, then \eqref{bound theta_k pre-ite} implies \begin{align*}
		\max_{\ell_{m}\le r\le k}\vert\theta_{r}\vert\le A_{m}+cA_{m}\sum_{\ell=\ell_{m}+1}^{k}\eta_{\ell} \prod_{j=\ell_{m}+1}^{\ell-1}(1+c\eta_{j})=A_{m}\prod_{\ell=\ell_{m}+1}^{k}(1+c\eta_{\ell}),
	\end{align*}
	where the inequality follows from induction hypothesis and the last step follows from a telescope sum. Actually, let $P_{k}\coloneqq \prod_{\ell=\ell_{m}+1}^{k}(1+c\eta_{\ell})$, then $P_{\ell_{m}}=1$, and $P_{k}=P_{k-1}(1+c\eta_{k})$ gives that $P_{k}-P_{k-1}=c\eta_{k}P_{k-1}$, summing over $\ell_{m}+1$ to $k$ yields \begin{align*}
		\prod_{\ell=\ell_{m}+1}^{k}(1+c\eta_{\ell})=1+c\sum_{\ell=\ell_{m}+1}^{k}\eta_{\ell}\prod_{j=\ell_{m}+1}^{\ell-1}(1+c\eta_{j}).
	\end{align*}
\end{proof}
\subsection{Proof of Lemma \ref{LemNum}}\label{SubsecStepsizes}

\begin{proof}[Proof of Lemma \ref{LemNum}]
	\underline{Proof of \eqref{LimNum1}}:  Direct computation gives \begin{align*}
		\frac{\sum_{1\le i<j\le n}\eta_{i}^{-3/4}\eta_{j}^{-1}}{\left(\sum_{k=1}^{n}\eta_{k}^{-1}\right)^{2}}\le\frac{\left(\sum_{i=1}^{n}\eta_{i}^{-3/4}\right)\left(\sum_{j=1}^{n}\eta_{j}^{-1}\right)}{\left(\sum_{k=1}^{n}\eta_{k}^{-1}\right)^{2}}=\frac{\sum_{k=1}^{n}\eta_{k}^{-3/4}}{\sum_{k=1}^{n}\eta_{k}^{-1}}.
	\end{align*}
	For arbitrarily  $\epsilon\in(0,1)$, monotonicity allows us to take some positive integer $n_{0}(\epsilon)=:n_{0}$, such that $\eta_{k}< \epsilon^{4}/16$ for all $k\ge n_{0}$. For $n>n_{0}$, splitting the sum by $n_{0}$ and dividing both numerator and denominator by $\sum_{k=n_{0}+1}^{n}\eta_{k}^{-1}$, we have
	\begin{align*}
		\frac{\sum_{k=1}^{n}\eta_{k}^{-3/4}}{\sum_{k=1}^{n}\eta_{k}^{-1}} \le \frac{\sum_{k=1}^{n_{0}}\eta_{k}^{-3/4}}{\sum_{k=n_{0}+1}^{n}\eta_{k}^{-1}}+\frac{\sum_{k=n_{0}+1}^{n}\eta_{k}^{-3/4}}{\sum_{k=n_{0}+1}^{n}\eta_{k}^{-1}}.
	\end{align*} 
	Since the first term is decreasing in $n$, enlarge $n_{0}$ to bound it by $\epsilon/2$. For the second term, the monotonicity yields $\sum_{k=n_{0}+1}^{n}\eta_{k}^{-1}\ge \eta_{n_{0}}^{-1/4}\sum_{k=n_{0}+1}^{n}\eta_{k}^{-3/4}$, which gives \begin{align*}
		\frac{\sum_{k=n_{0}+1}^{n}\eta_{k}^{-3/4}}{\sum_{k=n_{0}+1}^{n}\eta_{k}^{-1}}\le\eta_{n_{0}}^{1/4} \frac{\sum_{k=n_{0}+1}^{n}\eta_{k}^{-3/4}}{\sum_{k=n_{0}+1}^{n}\eta_{k}^{-3/4}}\le\frac{\epsilon}{2}.
	\end{align*}
	Combine estimations above, for all $n\ge n_{0}$, one has \begin{align*}
		\frac{\sum_{1\le i<j\le n}\eta_{i}^{-3/4}\eta_{j}^{-1}}{\left(\sum_{k=1}^{n}\eta_{k}^{-1}\right)^{2}}\le\frac{\sum_{k=1}^{n}\eta_{k}^{-3/4}}{\sum_{k=1}^{n}\eta_{k}^{-1}} <\epsilon.
	\end{align*}
	
	\underline{Proof of \eqref{bound SumEtaExp Num2}}: Denote $s_{n}\coloneqq\sum_{k=2}^{n}e^{-C_{1}(t_{n}-t_{k})}\eta_{k}^{1+q}$, then \eqref{bound SumEtaExp Num2} will follow from contraction if we can find some constants $c_{1}\in(0,1)$,  $c_{2}>0$, such that 
	\begin{align}\label{SumEtakWithNegExpReducedform}
		s_{n}\le c_{1}s_{n}+c_{2}\eta_{n}^{q},\quad\forall n\ge 2.
	\end{align}
	For proving \eqref{SumEtakWithNegExpReducedform}, note by monotonicity of $\{\eta_{k}\}_{k\ge 1}$ that \begin{align*}
		s_{n}\le e^{C_{1}\eta_{1}}\sum_{k=2}^{n}e^{-C_{1}(t_{n}-t_{k-1})}\eta_{k}^{1+q}\le C_{1}^{-1}e^{C_{1}\eta_{1}}\sum_{k=2}^{n}\eta_{k}^{q}[e^{-C_{1}(t_{n}-t_{k})}-e^{-C_{1}(t_{n}-t_{k-1})}],
	\end{align*}
	where the last inequality follows from the estimation\begin{align*}
		e^{-C_{1}(t_{n}-t_{n-1})}\eta_{k}^{1+q}\le e^{-C_{1}(t_{n}-t_{k-1})}\eta_{k}^{q}C_{1}^{-1}[e^{C_{1}\eta_{k}}-1].
	\end{align*}
	By adding and subtracting  $\eta_{k-1}^{q}e^{-C_{1}(t_{n}-t_{k-1})}$ in  summand, the sum becomes 
	\begin{align*}
		&\sum_{k=2}^{n}\eta_{k}^{q}[e^{-C_{1}(t_{n}-t_{k})}-e^{-C_{1}(t_{n}-t_{k-1})}]\\
		=&\sum_{k=2}^{n}[e^{-C_{1}(t_{n}-t_{k})}\eta_{k}^{q}-e^{-C_{1}(t_{n}-t_{k-1})}\eta_{k-1}^{q}]+\sum_{k=2}^{n}e^{-C_{1}(t_{n}-t_{k-1})}(\eta_{k-1}^{q}-\eta_{k}^{q})\\
		=&[\eta_{n}^{q}-e^{-C_{2}(t_{n}-t_{1})}\eta_{1}^{q}]+\sum_{k=2}^{n}e^{-C_{1}(t_{n}-t_{k-1})}(\eta_{k-1}^{q}-\eta_{k}^{q})\\
		\le&\eta_{n}^{q}+c\sum_{k=2}^{n}e^{-C_{1}(t_{n}-t_{k-1})}\eta_{k}^{1+q},
	\end{align*} 
	where the last inequality follows from the estimation\begin{align*}
		\eta_{k-1}^{q}-\eta_{k}^{q}\le \eta_{k}^{q-1}(\eta_{k-1}-\eta_{k})\le c\eta_{k}^{1+q}
	\end{align*}
	due to \eqref{con etak-diff ine} for $q\in[0,1]$. Since $e^{-C_{1}(t_{n}-t_{k-1})}=e^{-C_{1}(t_{n}-t_{k})}e^{-C_{1}\eta_{k}}\le e^{-C_{1}(t_{n}-t_{k})}$, one has \begin{align*}
		s_{n}\le C_{1}^{-1}e^{C_{1}\eta_{1}}\left[\eta_{n}^{q}+cs_{n}\right].
	\end{align*} 
	Note  that the constant $c$ in \eqref{con etak-diff ine}  can   be sufficiently small. Taking $c_{1}=C_{1}^{-1}e^{C_{1}\eta_{1}}c$, $c_{2}=C_{1}^{-1}e^{C_{1}\eta_{1}}$, choosing $c$ such that $c_{1}\in(0,1)$, we obtain \eqref{SumEtakWithNegExpReducedform}  as desired.

	\underline{Proof of \eqref{LimNum3}}: In the case $p=1$, consider $0<a<1$, by   H\"{o}lder's inequality for conjugate exponents $1/a$ and $1/(1-a)$, one has \begin{align*}
		\frac{\sum_{k=1}^{n}\eta_{k}^{1/2}}{\sqrt{T_{n}}}=\frac{\sum_{k=1}^{n}\eta_{k}^{a+1/2}\eta_{k}^{-a}}{\sqrt{T_{n}}}\le \frac{\left[\sum_{k=1}^{n}\eta_{k}^{(a+1/2)/(1-a)}\right]^{1-a}}{(T_{n})^{-a+1/2}}.
	\end{align*}
	To make the last term tend to $0$, constant $a$ should be satisfy\begin{align}\label{Ine}
		\left\{\begin{array}{l}
			-a+1/2>0,\\
			\frac{(2a+1)/2}{1-a}\ge2-\epsilon,
		\end{array}\right.
	\end{align}
	where the $\epsilon\in(0,1)$ is a given constant in Assumption \ref{Asp Etak-1}. Solving \eqref{Ine} to get
	\begin{align*}
		\frac{(3-2\epsilon)/2}{3-\epsilon}<a<\frac{1}{2}=\frac{(3-\epsilon)/2}{3-\epsilon}.
	\end{align*}
	Taking $a=(6-3\epsilon)/[4(3-\epsilon)]$ so that \eqref{Ine} applies, we see from \eqref{con etak-strongR-M} and \eqref{con etak-rate of etanTn} that \begin{align*}
		\frac{\sum_{k=1}^{n}\eta_{k}^{1/2}}{\sqrt{T_{n}}}\le \frac{\left[\sum_{k=1}^{n}\eta_{k}^{(a+1/2)/(1-a)}\right]^{1-a}}{(T_{n})^{-a+1/2}}\le \frac{C}{(T_{n})^{-a+1/2}}\to 0.
	\end{align*}
	For case $p>1$,  the monotonicity of $\{\eta_{k}\}_{k\ge1}$  yields $\sum_{k=1}^{n}\eta_{k}^{p/2}\le \eta_{1}^{(p-1)/2}\sum_{k=1}^{n}\eta_{k}^{1/2}$, and the result follows from that of case $p=1$.\\
	
	\underline{Proof of \eqref{LimNum4}}: The first limit follows immediately from monotonicity and \eqref{con etak-rate of etanTn}: \begin{align*}
		\frac{1}{T_{n}^{2}} \sum_{k=1}^{n}\frac{1}{\eta_{k}^{2}} \le \frac{1}{T^2_n}\frac{1}{\eta_{n}}\sum_{k=1}^n \frac{1}{\eta_{k}}=\frac{1}{\eta_{n} T_{n}} \rightarrow 0.
	\end{align*}  
	For the second limit,  note that the term $e^{-C(t_{j}-t_{i})}$ will not work when $t_{j}-t_{i}$ is small, we need to split the sum according to the size of $t_{j}-t_{i}$. Let $L_{n}\to \infty$ satisfy \eqref{con etak-rate of etanTn}.  Note by monotonicity of $\{\eta_{k}\}_{k\ge 1}$ that for all $i<j\le n$, $t_{j}-t_{i}\ge (j-i)\eta_{n}$. For each $i$, $j-i\ge \eta_{n}^{-1}L_{n}$ will imply $t_{j}-t_{i}\ge L_{n}$, hence for each $i$, $\sharp\{j:t_{j}-t_{i}< L_{n}\}\le \eta_{n}^{-1}L_{n}$, which gives \begin{align}\label{bound tj-ti<Ln}
		\sum_{\substack{1\le i<j\le n,\\
				t_{j}-t_{i}< L_{n}}}\frac{1}{\eta_{i}}\frac{1}{\eta_{j}}e^{-C(t_{j}-t_{i})}\le \eta_{n}^{-2}L_{n}\sum_{k=1}^{n}\eta_{k}^{-1}.
	\end{align}
	On the other hand, if $t_{j}-t_{i}\ge L_{n}$, then $e^{-C(t_{j}-t_{i})}\le e^{-CL_{n}}$, so \begin{align}\label{bound tj-ti>Ln}
		\sum_{\substack{1\le i<j\le n,\\
				t_{j}-t_{i}\ge L_{n}}}\frac{1}{\eta_{i}}\frac{1}{\eta_{j}}e^{-C(t_{j}-t_{i})}\le e^{-CL_{n}}\sum_{1\le i<j\le n}\eta_{i}^{-1}\eta_{j}^{-1},
	\end{align}
	and the latter sum is bounded by $T_{n}^{2}$. Using  \eqref{bound tj-ti<Ln} and \eqref{bound tj-ti>Ln} to get
	\begin{align*}
		\frac{1}{T_{n}^{2}}\sum_{1\le i<j\le n}\eta_{i}^{-1}\eta_{j}^{-1}e^{-C(t_{j}-t_{i})}\le C_{1}\frac{L_{n}}{\eta_{n}^{2}T_{n}}+C_{2}e^{-CL_{n}}\to 0
	\end{align*} 
	as $n\to\infty$, where \eqref{con etak-rate of etanTn} in Assumption \ref{Asp Etak-1} has been applied in the last step.  
\end{proof}
\subsection{Proof of Lemma \ref{LemRem}}\label{SubsecProofLemRem}
By direct computation, the remainder $R_{n}(m)$ in decomposition \eqref{def Decomposition} admits the decomposition \eqref{def R_n(m)} with \begin{align*}
	R_{n}^{(0)}(m)\coloneqq&\frac{1}{\sqrt{T_{n}}}\sum_{k=0}^{m-1}\frac{1}{\eta_{k+1}}[\varphi(\theta_{k+1})-\varphi(\theta_{k})],\\
	R_{n}^{(1)}(m) \coloneqq&\frac{1}{2\sqrt{T_{n}}}\sum_{k=0}^{m-1}\langle\nabla^{2}\varphi(\theta_{k}),\sigma\sigma^{\intercal}-(\sigma \xi_{k+1})(\sigma\xi_{k+1})^{\intercal}\rangle_{\text{HS}}\\
	&+\frac{1}{2\sqrt{T_{n}}}\sum_{k=0}^{m-1}\eta_{k+1}^{1/2}\langle\nabla^{2}\varphi(\theta_{k}),b(\theta_{k})(\sigma\xi_{k+1})^{\intercal}+(\sigma \xi_{k+1})b(\theta_{k})^{\intercal}\rangle_{\text{HS}},\\
	R_{n}^{(2)}(m)\coloneqq&\frac{1}{2\sqrt{T_{n}}}\sum_{k=0}^{m-1}\eta_{k+1}\langle\nabla^{2}\varphi(\theta_{k}),b(\theta_{k})b(\theta_{k})^{\intercal}\rangle_{\text{HS}},\\	
	R_{n}^{(3)}(m) \coloneqq&\frac{1}{\sqrt{T_{n}}}\sum_{k=0}^{m-1}\frac{1}{\eta_{k+1}}\mathcal{R}_{k+1}.
\end{align*}
\begin{proof} [Proof of Lemma \ref{LemRem}]
	We can actually prove $\sup_{t\in[0,1]}\vert R_{n}^{(i)}(\lfloor nt\rfloor)\vert\to 0$ for $i=1,2,3$ under  Assumptions \ref{Asp Coeff}, \ref{Asp Etak-0} and \ref{Asp Etak-1}; for $i=0$, we prove different results depending on whether Assumption \ref{Asp Etak-2} is imposed.  Motivated by Markov's inequality, we instead prove convergence in $L^{1}$ or $L^{2}$. We have nothing to do when   $t=0$, so we assume without loss of generality that $\lfloor nt\rfloor\ge 1$ below. 
	
	For $R_{n}^{(0)}(\lfloor nt\rfloor)$, the Abel transform (see for example \ocite{Rud64}*{Theorem 3.41}) yields 
	\begin{align}
		&\sum_{k=1}^{m}\frac{1}{\eta_{k}}[\varphi(\theta_{k})-\varphi(\theta_{k-1})]\nonumber\\
		=&\sum_{k=1}^{m-1}\left(\frac{1}{\eta_{k}}-\frac{1}{\eta_{k+1}}\right)\left(\sum_{i=1}^{k}[\varphi(\theta_{i})-\varphi(\theta_{i-1})]\right)+\frac{1}{\eta_{m}}\left(\sum_{i=1}^{m}[\varphi(\theta_{i})-\varphi(\theta_{i-1})]\right)\nonumber\\
		=&\sum_{k=1}^{m-1}\left(\frac{1}{\eta_{k}}-\frac{1}{\eta_{k+1}}\right)[\varphi(\theta_{k})-\varphi(\theta_{0})]+\frac{1}{\eta_{m}}[\varphi(\theta_{m})-\varphi(\theta_{0})]\label{abel transform R0}
	\end{align} 
	for $1\le m\le n$. Taking $m=\lfloor nt\rfloor$ yields
	\begin{align*}
		&\sup_{t\in[0,1]}\left\vert \sum_{k=1}^{\lfloor nt\rfloor}\frac{1}{\eta_{k}}[\varphi(\theta_{k})-\varphi(\theta_{k-1})]\right\vert\\
		\le&\sup_{t\in[0,1]}\sum_{k=1}^{\lfloor nt\rfloor-1}\left\vert\left(\frac{1}{\eta_{k}}-\frac{1}{\eta_{k+1}}\right)[\varphi(\theta_{k})-\varphi(\theta_{0})]\right\vert+\frac{1}{\eta_{n}}\sup_{t\in[0,1]}\left\vert\varphi(\theta_{\lfloor nt\rfloor})-\varphi(\theta_{0})\right\vert\\
		\le&\sum_{k=1}^{n-1}\left(\frac{1}{\eta_{k+1}}-\frac{1}{\eta_{k}}\right) \vert \varphi(\theta_{k})-\varphi(\theta_{0})\vert+\frac{1}{\eta_{n}}\sup_{t\in[0,1]}\left\vert\varphi(\theta_{\lfloor nt\rfloor})\right\vert +\frac{1}{\eta_{n}}\left\vert\varphi(\theta_{0})\right\vert. 
	\end{align*}
	By \eqref{bound |grad_k phi(x)|} and \eqref{bound E|b(theta_k)|^p}, we have 
	\begin{align}\label{bound E|phi(theta_k)-phi(theta_0)|}
		\mathbb{E} \vert\varphi(\theta_{k})-\varphi(\theta_{0})\vert\le  C\mathbb{E} (1+\vert\theta_{k}\vert^{2}+\vert x\vert^{2}) \le C(1+\vert x\vert^2)  ,
	\end{align} 
	for all $k\ge 1$ and for some constant $C$  independent of $k$. Since $\sum_{k=1}^{n-1}\left(\frac{1}{\eta_{k+1}}-\frac{1}{\eta_{k}}\right)\le\frac{1}{\eta_{n}}$, \eqref{bound E|phi(theta_k)-phi(theta_0)|} implies that \begin{align}\label{bound sum(eta_k+1-eta_k)E|phi-phi|}
		\sum_{k=1}^{n-1}\left(\frac{1}{\eta_{k+1}}-\frac{1}{\eta_{k}}\right) \mathbb{E}\vert \varphi(\theta_{k})-\varphi(\theta_{0})\vert\le \frac{C}{\eta_{n}}(1+\vert x\vert^{2}).
	\end{align}
	Meanwhile, by H\"{o}lder's inequality, \eqref{bound |grad_k phi(x)|} and \eqref{bound theta_k^2p uniformly} imply that \begin{align*}
		\mathbb{E}\left[\sup_{t\in[0,1]}\left\vert\varphi(\theta_{\lfloor nt\rfloor})\right\vert\right]\le& C\mathbb{E}\left[\max_{0\le k\le n}\left(1+\vert \theta_{k}\vert^{2} \right)\right]\\
		\le&\mathbb{E}\left[\max_{0\le k\le n}\left(1+\vert \theta_{k}\vert^{2}\right)^{p}\right]^{1/p}\le C(1+\vert x\vert^{2})(1+t_{n})^{1/p}
	\end{align*}
	for $p\ge 1$ as in Assumption \ref{Asp Etak-2}. Combine arguments above, by  \eqref{con etak-rate uniform}, as $n\to\infty$, \begin{align*}
		\mathbb{E}\left[\frac{1}{\sqrt{T_{n}}}\sup_{t\in[0,1]}\left\vert \sum_{k=1}^{\lfloor nt\rfloor}\frac{1}{\eta_{k}}[\varphi(\theta_{k})-\varphi(\theta_{k-1})]\right\vert\right]\le C(1+\vert x\vert^{2})\left[\frac{(1+t_{n})^{1/p}}{\eta_{n}\sqrt{T_{n}}}\right]\to 0.
	\end{align*} 
	
	For $R_{n}^{(0)}(n)$,  Assumption \ref{Asp Etak-2} is not needed. Taking $m=n$ in \eqref{abel transform R0} yields
	\begin{align*}
		\left\vert\sum_{k=0}^{n-1}\frac{1}{\eta_{k+1}}[\varphi(\theta_{k+1})-\varphi(\theta_{k})]\right\vert\le \sum_{k=1}^{n-1}\left(\frac{1}{\eta_{k+1}}-\frac{1}{\eta_{k}}\right) \vert \varphi(\theta_{k})-\varphi(\theta_{0})\vert+\frac{1}{\eta_{n}}\vert\varphi(\theta_{n})-\varphi(\theta_{0})\vert,
	\end{align*}
	by \eqref{bound E|phi(theta_k)-phi(theta_0)|} and \eqref{bound sum(eta_k+1-eta_k)E|phi-phi|}, we have
	\begin{align*} 
		\mathbb{E}\vert R_{n}^{(0)}(n)\vert\le\frac{1}{\sqrt{T_{n}}}\mathbb{E}\left\vert\sum_{k=0}^{n-1}\frac{1}{\eta_{k+1}}[\varphi(\theta_{k+1})-\varphi(\theta_{k})]\right\vert\le \frac{C(1+\vert x \vert^{2})}{\eta_{n}\sqrt{T_{n}}}\to0,
	\end{align*}
	where the  last step is by \eqref{con etak-rate of etanTn} in  Assumption \ref{Asp Etak-1}.
	
	For $R_{n}^{(1)}(\lfloor nt\rfloor)$, note that $\{R_{n}^{(1)}(m)\}_{m\ge 1}$ is a martingale. By Doob's $L^{2}$ inequality \ocite{RY99}*{Theorem \uppercase\expandafter{\romannumeral2\relax}.(1.7)}, one has\begin{align*}
		\mathbb{E}\left[\sup_{t\in[0,1]}\vert R_{n}^{(1)}(\lfloor nt\rfloor)\vert^{2}\right]\le \mathbb{E}\left[\sup_{1\le m\le n}\vert R_{n}^{(1)}(m)\vert^{2}\right]\le C\mathbb{E}\vert R_{n}^{(1)}(n)\vert^{2}.
	\end{align*}
	Since $\mathbb{E}\vert R_{n}^{(1)}(n)\vert^{2}\le 2\mathbb{E}\vert R_{n}^{(1),1}(n)\vert^{2}+2\mathbb{E}\vert R_{n}^{(1),2}(n)\vert^{2}$, where   \begin{align*}
		R_{n}^{(1),1}(n)&\coloneqq\frac{1}{2\sqrt{T_{n}}}\sum_{k=0}^{n-1}\langle\nabla^{2}\varphi(\theta_{k}),\sigma\sigma^{\intercal}-(\sigma \xi_{k+1})(\sigma\xi_{k+1})^{\intercal}\rangle_{\text{HS}},\\
		R_{n}^{(1),2}(n)&\coloneqq\frac{1}{2\sqrt{T_{n}}}\sum_{k=0}^{n-1}\eta_{k+1}^{1/2}\langle\nabla^{2}\varphi(\theta_{k}),b(\theta_{k})(\sigma\xi_{k+1})^{\intercal}+(\sigma \xi_{k+1})b(\theta_{k})^{\intercal}\rangle_{\text{HS}},
	\end{align*}
	we turn to bound the latter two expectations.  Note that \begin{align*}
		\mathbb{E}\left[\sum_{k=0}^{n-1}\langle\nabla^{2}\varphi(\theta_{k}),\sigma\sigma^{\intercal}-(\sigma \xi_{k+1})(\sigma\xi_{k+1})^{\intercal}\rangle_{\text{HS}}\right]^{2}=\mathbb{E}\sum_{k=0}^{n-1}\langle\nabla^{2}\varphi(\theta_{k}),\sigma\sigma^{\intercal}-(\sigma\xi_{k+1})(\sigma\xi_{k+1})^{\intercal}\rangle_{\text{HS}}^{2},
	\end{align*}
	where the cross term in expansion \begin{align*}
		\mathbb{E}\sum_{\substack{
				0\le i<j\le n-1
		}}\langle\nabla^{2}\varphi(\theta_{i}),\sigma\sigma^{\intercal}-(\sigma\xi_{i+1})(\sigma\xi_{i+1})^{\intercal}\rangle_{\text{HS}}\langle\nabla^{2}\varphi(\theta_{j}),\sigma\sigma^{\intercal}-(\sigma\xi_{j+1})(\sigma\xi_{j+1})^{\intercal}\rangle_{\text{HS}}=0
	\end{align*}
	by conditioning on $\theta_{i}$ and using the fact $\mathbb{E}[ \sigma\sigma^{\intercal}-(\sigma\xi_{k+1})(\sigma\xi_{k+1})^{\intercal}] =0$. Therefore,  \begin{align*}
		\mathbb{E}\vert R_{n}^{(1),1}(n)\vert^{2}\le &\frac{1}{4T_{n}}\mathbb{E}\sum_{k=0}^{n-1}\langle\nabla^{2}\varphi(\theta_{k}),\sigma\sigma^{\intercal}-(\sigma\xi_{k+1})(\sigma\xi_{k+1})^{\intercal}\rangle_{\text{HS}}^{2}\\
		\le&\frac{1}{4T_{n}}\sum_{k=0}^{n-1}\mathbb{E}\vert \nabla^{2}\varphi(\theta_{k})\vert^{2}\mathbb{E}\vert \sigma\sigma^{\intercal}-(\sigma\xi_{k+1})(\sigma\xi_{k+1})^{\intercal}\vert^{2}\\
		\le& \frac{Cn}{T_{n}}(1+\vert x\vert^{8}),
	\end{align*}
	where the third inequality follows from \eqref{bound |grad_k phi(x)|} and moment estimation. Similarly, one has \begin{align*}
		\mathbb{E}\vert R_{n}^{(1),2}(n)\vert^{2}=\frac{1}{4T_{n}}\sum_{k=0}^{n-1}\eta_{k+1}\mathbb{E}[\langle\nabla^{2}\varphi(\theta_{k}),b(\theta_{k})(\sigma\xi_{k+1})^{\intercal}+(\sigma \xi_{k+1})b(\theta_{k})^{\intercal}\rangle_{\text{HS}}]^{2},
	\end{align*}
	where the cross terms vanished by conditioning. By \eqref{bound |grad_k phi(x)|} and moment estimation, we have 
	\begin{align*}
		\mathbb{E}[\langle\nabla^{2}\varphi(\theta_{k}),b(\theta_{k})(\sigma\xi_{k+1})^{\intercal}+(\sigma \xi_{k+1})b(\theta_{k})^{\intercal}\rangle_{\text{HS}}]^{2}\le C(1+\vert x\vert^{10}),
	\end{align*}
	which gives 
	\begin{align*} 
		\mathbb{E}\vert R_{n}^{(1),2}(n)\vert^{2}\le C(1+\vert x\vert^{10})\frac{1}{T_{n}}\sum_{k=1}^{n}\eta_{k}.
	\end{align*} 
	Combining arguments above, we arrive at \begin{align*}
		\mathbb{E}\sup_{t\in[0,1]}\vert R_{n}^{(1)}(\lfloor nt\rfloor)\vert^{2}\le& C(1+\vert x\vert^{10})\left(\frac{n}{T_{n}}+\frac{1}{T_{n}}\sum_{k=1}^{n}\eta_{k}\right)\\
		\le& C\left(\frac{\sqrt{\sum_{k=1}^{n}\eta_{k}^{2}}}{\sqrt{n}}+\frac{\sum_{k=1}^{n}\eta_{k}}{T_{n}}\right)\to 0
	\end{align*}
	where the first inequality follows from  QM-HM inequality  \ocite{SS18}*{(3.5)}, which says that  $\frac{n}{\sum_{k=1}^{n}a_{k}^{-1}}\le\left(\frac{\sum_{k=1}^{n}a_{k}^{2}}{n}\right)^{1/2}$ for positive sequence $\{a_{k}\}_{k=1}^{n}$, and the last step is because  the finiteness of $\sum_{k=1}^{\infty}\eta_{k}^{2}$ and Lemma \ref{LemNum}.
	
	For $R_{n}^{(2)}(\lfloor nt\rfloor)$, note that \begin{align*}
		\mathbb{E}\sup_{t\in[0,1]}\vert R_{n}^{(2)}(\lfloor nt\rfloor)\vert\le& \frac{1}{2\sqrt{T_{n}}}\mathbb{E}\left[\sup_{t\in[0,1]}\sum_{k=0}^{\lfloor nt\rfloor-1}\eta_{k+1}\left\vert \langle\nabla^{2}\varphi(\theta_{k}),b(\theta_{k})b(\theta_{k})^{\intercal}\rangle_{\text{HS}}\right\vert\right]\\
		\le &\frac{C}{\sqrt{T_{n}}}\sum_{k=0}^{n-1}\eta_{k+1}\mathbb{E}\left[\left\vert  \nabla^{2}\varphi(\theta_{k})\right\vert \left\vert b(\theta_{k})\right\vert^{2}\right].
	\end{align*}
	The growth conditions  \eqref{ine LinearGrowth} and \eqref{bound |grad_k phi(x)|} yield \begin{align*}
		\mathbb{E}\left[\left\vert  \nabla^{2}\varphi(\theta_{k})\right\vert \left\vert b(\theta_{k})\right\vert^{2}\right]\le C(1+\mathbb{E}\vert \theta_{k}\vert ^{6}) \le C(1+\vert x\vert ^{6}),
	\end{align*}
	which, combining \eqref{LimNum3}, gives\begin{align*}
		\mathbb{E}\sup_{t\in[0,1]}\vert R_{n}^{(2)}(\lfloor nt\rfloor)\vert\le C(1+\vert x\vert^{6})\frac{1}{\sqrt{T_{n}}}\sum_{k=1}^{n}\eta_{k}\to 0.
	\end{align*}
 
	For $R_{n}^{(3)}(\lfloor nt\rfloor)$, by a similar argument, one has \begin{align*}
		\mathbb{E}\sup_{t\in[0,1]}\vert R_{n}^{(3)}(\lfloor nt\rfloor)\vert\le \frac{1}{\sqrt{T_{n}}}\sum_{k=0}^{n-1}\frac{1}{\eta_{k+1}}\mathbb{E}\vert \mathcal{R}_{k+1}\vert.
	\end{align*}
	We use \begin{align*}
		\mathbb{E}\vert \mathcal{R}_{k+1}\vert\le \mathbb{E}\vert\mathcal{R}_{k+1}1_{\{\vert\Delta \theta_{k+1}\vert\le 1\}}\vert+\mathbb{E}\vert \mathcal{R}_{k+1}1_{\{\vert\Delta \theta_{k+1}\vert> 1\}}\vert,
	\end{align*} where $\Delta\theta_{k+1}=\theta_{k+1}-\theta_{k}$ satisfies \begin{align*}
		\mathbb{E}\vert \Delta\theta_{k+1}\vert\le \eta_{k+1}\mathbb{E}\vert b(\theta_{k})\vert+\sqrt{\eta_{k+1}}\mathbb{E}\vert\sigma\xi_{k+1}\vert \le C(1+\vert x\vert)\eta_{k+1}^{1/2}.
	\end{align*} 
	On the 
	set $\{\vert\Delta \theta_{k+1}\vert\le 1\}$, we use \eqref{bound DiffHesPhi} and H\"{o}lder's inequality to get
	\begin{align*}
		&\mathbb{E}\left\vert\mathcal{R}_{k+1}1_{\{\vert \Delta\theta_{k+1}\vert\le 1\}}\right\vert\\
		\le&\int_{0}^{1}\int_{0}^{1}r^{2}s\mathbb{E}\frac{\vert \nabla^{2}\varphi(\theta_{k}+r  s (\Delta\theta_{k+1}))-\nabla^{2}\varphi(\theta_{k})\vert}{\vert rs \Delta\theta_{k+1}\vert}\vert\Delta\theta_{k+1}\vert^{3}1_{\{\vert \Delta\theta_{k+1}\vert\le 1\}}{\rm d}r{\rm d}s\\
		\le&C\mathbb{E}[(1+\vert \theta_{k}+\Delta\theta_{k+1}\vert^{5})\vert\Delta\theta_{k+1}\vert^{3}1_{\{\vert \Delta\theta_{k+1}\vert\le 1\}}]\\
		\le&C\left[\mathbb{E} (1+\vert \theta_{k+1}\vert^{5})^{2}1_{\{\vert \Delta\theta_{k+1}\vert\le 1\}}\right]^{1/2}[\mathbb{E} \vert \Delta\theta_{k+1}\vert^{6}]^{1/2}\\
		\le&C(1+\vert x\vert^{8})\eta_{k+1}^{3/2},
	\end{align*} 
	while on the set $\{\vert\Delta \theta_{k+1}\vert>1\}$, 
	\begin{align*}
		&\mathbb{E}\left\vert\mathcal{R}_{k+1}1_{\{\vert \Delta\theta_{k+1}\vert> 1\}}\right\vert\\
		\le&\int_{0}^{1}\int_{0}^{1}r\mathbb{E}\left[(\vert \nabla^{2}\varphi(\theta_{k}+rs\Delta\theta_{k+1})\vert+\vert \nabla^{2}\varphi(\theta_{k}))\vert)\vert\Delta\theta_{k+1}\vert^{2}1_{\{\vert \Delta\theta_{k+1}\vert> 1\}}\right]{\rm d}r{\rm d}s\\
		\le&C\mathbb{E}\left[\left(1+\vert \theta_{k}\vert^{4}+\vert\theta_{k}+\Delta\theta_{k+1}\vert^{4}\right)\vert\Delta\theta_{k+1}\vert^{2}1_{\{\vert \Delta\theta_{k+1}\vert> 1\}}\right]\\
		\le&C\left[\mathbb{E}\left(1+\vert \theta_{k}\vert^{4}+\vert\theta_{k}+\Delta\theta_{k+1}\vert^{4}\right)^{2}\vert\Delta\theta_{k+1}\vert^{4}\right]^{1/2}\left[\mathbb{E}1_{\{\vert \Delta\theta_{k+1}\vert> 1\}}\right]^{1/2}\\
		\le&C(1+\vert x\vert^{12})\eta_{k+1}^{2},
	\end{align*} 
	where the last step used the following estimation by Chebyshev and \eqref{bound E|b(theta_k)|^p}:\begin{align*}
		\mathbb{P}\left(\vert \Delta\theta_{k+1}\vert>1\right)=\mathbb{P}\left(\vert\Delta\theta_{k+1}\vert^{4}>1\right)\le \mathbb{E}\vert\Delta\theta_{k+1}\vert^{4}\le C(1+\vert x\vert^{4})\eta_{k+1}^{2}.
	\end{align*} 
	Combine estimations above and applying Lemma \ref{LemNum}, we arrive at
	\begin{align*}
		\mathbb{E}\sup_{t\in[0,1]}\vert R_{n}^{(3)}(\lfloor nt\rfloor)\vert\le & \frac{1}{\sqrt{T_{n}}}\sum_{k=0}^{n-1}\frac{1}{\eta_{k+1}}[\mathbb{E}\left\vert\mathcal{R}_{k+1}1_{\{\vert \Delta\theta_{k+1}\vert\le 1\}}\right\vert+\mathbb{E}\left\vert\mathcal{R}_{k+1}1_{\{\vert \Delta\theta_{k+1}\vert> 1\}}\right\vert]\\
		\le&\frac{1}{\sqrt{T_{n}}}\sum_{k=0}^{n-1}[C_{1}\eta_{k+1}^{1/2}+C_{2}\eta_{k+1}]\to 0.
	\end{align*}  
\end{proof}
  
\section{Proof of lemmas in Sections \ref{SectionProofCLT} and \ref{SectionProofW2distance}}\label{SectionProofLemmaThms}
\subsection{Proof of Lemma \ref{LemmaWeakComparison}}
\begin{proof}[Proof of Lemma \ref{LemmaWeakComparison}]
	Let $\{\widehat{X}_{t}\}_{t\ge t_{i}}$ be the solution to SDE \begin{align}\label{SDEstartfromTheta}
		{\rm d}\widehat{X}_{t}=b(\widehat{X}_{t}){\rm d}t+\sigma{\rm d}B_{t},\  \widehat{X}_{t_{i}}=\theta_{i}
	\end{align}
	where $b$ and $\sigma$ satisfies Assumption \ref{Asp Coeff}, $\{B_{t}\}_{t\ge0}$ is $d$-dimensional standard Brownian motion independent of $\theta_{i}$, and $\theta_{i}$ is given by \eqref{def EMDis}. By inserting the term $\mathbb{E}_{i}\vert \sigma^{\intercal}\nabla\varphi(\widehat{X}_{t_{j-1}})\vert^{2}$, triangle inequality gives 
	\begin{align}
		\label{SplitEphi-pi}\begin{split}
			&\left\vert\mathbb{E}_{i}[\vert\sigma^{\intercal}\nabla\varphi(\theta_{j-1})\vert^{2}-\pi(\vert\sigma^{\intercal}\nabla\varphi\vert^{2})]\right\vert\\
			\le&\vert \mathbb{E}_{i}\vert\sigma^{\intercal}\nabla\varphi(\theta_{j-1})\vert^{2}-\mathbb{E}_{i}\vert\sigma^{\intercal}\nabla\varphi(\widehat{X}_{t_{j-1}})\vert^{2}\vert+\vert \mathbb{E}_{i}\vert\sigma^{\intercal}\nabla\varphi(\widehat{X}_{t_{j-1}})\vert^{2}-\pi(\vert\sigma^{\intercal}\nabla\varphi\vert^{2})\vert.
		\end{split}
	\end{align} 
	
	For bounding the latter term in \eqref{SplitEphi-pi}, we use ergodicity of SDE \eqref{def SDE}. It is easy to see that $\widehat{X}_{t_{j-1}}\stackrel{d}{=}X_{t_{j-1}-t_{i}}$ with $X_{0}=\theta_{i}$ be given. The irreducibility of Markov process $\{X_{t}\}_{t\ge0}$ has been checked in \ocite{LTX22}*{Lemma 2.3}. Take Lyapunov function $V_{6}(x)=1+\vert x\vert^{6}$, the dissipativity \eqref{ine DispOneSide} gives for some constant $C,C_{1}>0$ that\begin{align*}
		\mathcal{A}V_{6}(x)=6\vert x\vert^{4}\langle x,b(x)\rangle+C\vert x\vert^{4}\le-3K_{1}\vert x\vert^{6}+C_{1}\vert x\vert^{4}\le -\frac{3}{2}K_{1}V_{6}(x)+q_{1}1_{A}(x),
	\end{align*}
	where the constant $K_{1}>0$ is given in \eqref{ine DispOneSide}, $q_{1}=C_{1}\left(1+(2C_{1})/(3K_{1})\right)^{2}+(3/2)K_{1}$, and $A=\{\vert x\vert\le\sqrt{1+(2C_{1})/(3K_{1})}\}$. Note that $\vert \nabla\varphi(x)\vert^{2}\le  CV_{6}(x)$, by \ocite{MT93}*{Theorem 6.1}, for arbitrary $x\in\mathbb{R}^{d}$ one has
	\begin{align*}
		\vert \mathbb{E}\vert \sigma^{\intercal}\nabla\varphi(X_{t})\vert^{2}-\pi(\vert\sigma^{\intercal}\nabla\varphi\vert^{2})\vert\le C_{1}V_{6}(x)e^{-C_{0}t}
	\end{align*}
	for some constant $C_{0}>0$, which implies that
	\begin{align*}
		\vert \mathbb{E}_{i}\vert\sigma^{\intercal}\nabla\varphi(\widehat{X}_{t_{j-1}})\vert^{2}-\pi(\vert\sigma^{\intercal}\nabla\varphi\vert^{2})\vert\le C_{1}V_{6}(\theta_{i})e^{-C_{0}(t_{j-1}-t_{i})}.
	\end{align*}

	For bounding another term in \eqref{SplitEphi-pi}, note by \eqref{bound |grad_k phi(x)|} that $\vert \nabla\vert\sigma^{\intercal}\nabla\varphi(x)\vert^{2}\vert\le C(1+\vert x\vert^{7})$, Taylor expansion and conditional  H\"{o}lder's inequality (see for example \ocite{Dur19}*{Page187})  imply
	\begin{align*}
		&\left\vert\mathbb{E}\left[\vert\sigma^{\intercal}\nabla\varphi(\theta_{j-1})\vert^{2}-\vert\sigma^{\intercal}\nabla\varphi(\widehat{X}_{t_{j-1}})\vert^{2}\mid \theta_{i}\right]\right\vert\\
		=&\left\vert\int_{0}^{1}\mathbb{E}\left[\nabla\vert\sigma^{\intercal}\nabla\varphi(\widehat{X}_{t_{j-1}}+r(\theta_{j-1}-\widehat{X}_{t_{j-1}}))\vert^{2}\cdot(\theta_{j-1}-\widehat{X}_{t_{j-1}}) \mid\theta_{i}\right]{\rm d}r\right\vert\\
		\le&C\int_{0}^{1}\left[\mathbb{E}[1+\vert \widehat{X}_{t_{j-1}}+r(\theta_{j-1}-\widehat{X}_{t_{j-1}})\vert^{14}\vert \theta_{i}]\right]^{1/2}{\rm d}r\left[\mathbb{E}[\vert\theta_{j-1}-\widehat{X}_{t_{j-1}}\vert^{2}\vert\theta_{i}]\right]^{1/2}\\
		\le&C(1+\vert\theta_{i}\vert^{7})\left[\mathbb{E}[\vert\theta_{j-1}-\widehat{X}_{t_{j-1}}\vert^{2}\vert\theta_{i}]\right]^{1/2},
	\end{align*}
	where the expectation $\mathbb{E}$ above is taken with respect to any coupling realization of $\theta_{j-1}$ and $\widehat{X}_{t_{j-1}}$, and the last inequality follows from \eqref{bound E|b(theta_k)|^p}. By taking a coupling attaining $\mathbf{W}_{2}$--distance and using the upper bound in Theorem \ref{TheoremW2distance}  with given $\theta_{i}$, it follows that \begin{align*}
		\vert \mathbb{E}_{i}\vert\sigma^{\intercal}\nabla\varphi(\theta_{j-1})\vert^{2}-\mathbb{E}_{i}\vert\sigma^{\intercal}\nabla\varphi(\widehat{X}_{t_{j-1}})\vert^{2}\vert\le C_{2}(1+\vert \theta_{i}\vert^{15/2})\eta_{j-1}^{1/4}.
	\end{align*}
\end{proof} 
\subsection{Proof of Lemma \ref{LemmaDifferenceCoupling}}
\begin{proof}[Proof of Lemma \ref{LemmaDifferenceCoupling}]
	Note that \begin{align*}
		{\rm d}Z_{t}=[b(\widetilde{X}_{t})-b(\widetilde{Y}_{t_{n-1}})]{\rm d}t+2\phi_{1}^{\delta}(Z_{t})\frac{z_{t}(\sigma^{-1}z_{t})^{\intercal}}{\vert \sigma^{-1}z_{t}\vert^{2}}{\rm d}B_{t}^{(1)},
	\end{align*}
	It\^{o}'s formula gives \begin{align*}
		{\rm d}\vert Z_{t}\vert^{2} =2\langle b(\widetilde{X}_{t})-b(\widetilde{Y}_{t_{n-1}}),Z_{t}\rangle{\rm d}t+4\phi_{1}^{\delta}(Z_{t})\frac{\vert Z_{t}\vert(\sigma^{-1}z_{t})^{\intercal}}{\vert \sigma^{-1}z_{t}\vert^{2}}{\rm d}B_{t}^{(1)}+4[\phi_{1}^{\delta}(Z_{t})]^{2}\frac{1}{\vert\sigma^{-1}z_{t}\vert^{2}}{\rm d}t.
	\end{align*}
	Take $\psi_{a}(r)\coloneqq(a+r)^{1/2}$, then $\psi_{a}\in C^{2}(\mathbb{R}_{+};\mathbb{R}_{+})$, $\psi_{a}^{\prime}(r)=\frac{1}{2(a+r)^{1/2}}$, $\psi_{a}^{\prime\prime}(r)=-\frac{1}{4(a+r)^{3/2}}$. Using It\^{o}'s formula again, we obtain
	\begin{align}\label{equation dPsi_a}
		\begin{split}
			d\psi_{a}(\vert Z_{t}\vert^{2})=&2\psi_{a}^{\prime}(\vert Z_{t}\vert^{2})\left[\langle b(\widetilde{X}_{t})-b(\widetilde{Y}_{t_{n-1}}),Z_{t}\rangle {\rm d}t+2\phi_{1}^{\delta}(Z_{t})\frac{\vert Z_{t}\vert(\sigma^{-1}z_{t})^{\intercal}}{\vert \sigma^{-1}z_{t}\vert^{2}}{\rm d}B_{t}^{(1)}\right] \\
			&+4\frac{[\phi_{1}^{\delta}(Z_{t})]^{2}}{\vert\sigma^{-1}z_{t}\vert^{2}}\left[\psi_{a}^{\prime}(\vert Z_{t}\vert^{2})+2\psi_{a}^{\prime\prime}(\vert Z_{t}\vert^{2}) \vert Z_{t}\vert^{2}\right]{\rm d}t.
		\end{split}
	\end{align} 
	For any $T\in[0,\eta_{n}]$, since $2r\psi_{a}^{\prime}(r^{2})=\frac{r}{(a+r^{2})^{1/2}}$, $\vert 2r\psi_{a}^{\prime}(r^{2})\vert\le 1$,  by (stochastic) dominated convergence theorem, one has
	\begin{align*}
		\lim_{a\to0}\int_{0}^{T}2\psi_{a}^{\prime}(\vert Z_{t}\vert^{2})\langle b(\widetilde{X}_{t})-b(\widetilde{Y}_{t_{n-1}}),Z_{t}\rangle {\rm d}t=&\int_{0}^{T}\langle b(\widetilde{X}_{t})-b(\widetilde{Y}_{t_{n-1}}),z_{t}\rangle {\rm d}t,\\
		\lim_{a\to0}\int_{0}^{T}4\psi_{a}^{\prime}(\vert Z_{t}\vert^{2})\phi_{1}^{\delta}(Z_{t})\frac{\vert Z_{t}\vert(\sigma^{-1}z_{t})^{\intercal}}{\vert \sigma^{-1}z_{t}\vert^{2}}{\rm d}B_{t}^{(1)}=&\int_{0}^{T}2\phi_{1}^{\delta}(Z_{t})\frac{ (\sigma^{-1}z_{t})^{\intercal}}{\vert \sigma^{-1}z_{t}\vert^{2}}{\rm d}B_{t}^{(1)};
	\end{align*}
	since $\psi_{a}^{\prime}(r^{2})+2r^{2}\psi_{a}^{\prime\prime}(r^{2})=\frac{a}{2(r^{2}+a)^{3/2}}\le\frac{a}{2r^{3}}$, $\phi_{1}^{\delta}(z)=0$ for $\vert z\vert<\frac{\delta}{2}$, it follows that \begin{align*}
		\lim_{a\to0}\int_{0}^{T}[\phi_{1}^{\delta}(Z_{t})]^{2}\left[\psi_{a}^{\prime}(\vert Z_{t}\vert^{2})+2\psi_{a}^{\prime\prime}(\vert Z_{t}\vert^{2})\vert Z_{t}\vert^{2}\right]\frac{1}{\vert\sigma^{-1}z_{t}\vert^{2}}{\rm d}t=0.
	\end{align*}
	The result then follows from letting $a\to0$ in \eqref{equation dPsi_a}. 
\end{proof}
\subsection{Proof of Lemma \ref{LemmaContraction}}
\begin{proof}[Proof of Lemma \ref{LemmaContraction}]
	We prove \eqref{ddtErho1} by bounding the right-hand side of inequality \eqref{dErho(X,Y)}. 
	
	\underline{Bounding $\frac{\rm d}{\rm dt}\mathbb{E} [\varepsilon V(X_{t})+\varepsilon V(Y_{t})]$: }
	By It\^{o}'s formula and Lyapunov condition, one has \begin{align}\label{dE|V(x)+V(y)|}
		{\rm d}\mathbb{E}[\varepsilon V(X_{t})+\varepsilon V(Y_{t})]\le -\lambda \mathbb{E}[\varepsilon V(X_{t})+\varepsilon V(Y_{t})]{\rm d}t+ q\varepsilon {\rm d}t,
	\end{align}
	where the constant $q$ is given by \eqref{ine Lyapnov} and the constant $\varepsilon$ will be selected later. 
	
	\underline{Bounding $\frac{\rm d}{\rm dt}\mathbb{E} f(\vert Z_{t}\vert)$: } Using It\^{o}'s formula again, it follows a.s. that \begin{align*}
		{\rm d}f(\vert Z_{t}\vert)= f^{\prime}(\vert Z_{t}\vert)\langle z_{t},b(\widetilde{X}_{t})-b(\widetilde{Y}_{t_{n-1}})\rangle {\rm d}t+2[\phi_{1}^{\delta}(Z_{t})]^{2}f^{\prime\prime}(\vert Z_{t}\vert )\frac{1}{\vert\sigma^{-1}z_{t}\vert^{2}}{\rm d}t+{\rm d}M_{t},
	\end{align*}
	where $\{M_{t}\}_{t\ge0}$ is a martingale. It implies that
	\begin{align}
		\label{bound dEf(Z)}  
		\frac{\rm d}{\rm dt}\mathbb{E} f(\vert Z_{t}\vert)=E_{1}+E_{2},
	\end{align}	where 
	\begin{align*}
		E_{1}&\coloneqq\mathbb{E}\left[ f^{\prime}(\vert Z_{t}\vert)\langle z_{t},b(\widetilde{Y}_{t})-b(\widetilde{Y}_{t_{n-1}})\rangle\right],\\
		E_{2}&\coloneqq\mathbb{E}\left[f^{\prime}(\vert Z_{t}\vert)\langle z_{t},b(\widetilde{X}_{t})-b(\widetilde{Y}_{t}) \rangle+2[\phi_{1}^{\delta}(Z_{t})]^{2}f^{\prime\prime}(\vert Z_{t}\vert)\frac{1}{\vert\sigma^{-1}z_{t}\vert^{2}} \right].
	\end{align*}
	For $E_{1}$, the Cauchy--Schwartz inequality, \eqref{boundedfprime} and Assumption \ref{Asp Coeff} give that \begin{align}\label{bound dEf(Z) Pt1}
		E_{1}\le L[\mathbb{E}\vert \widetilde{Y}_{t}-\widetilde{Y}_{t_{n-1}}\vert^{2}]^{1/2}\le C(1+\vert z\vert)\eta_{n}^{1/2},
	\end{align}
	where the last inequality follows from \eqref{WithinOneStep} and \eqref{bound E|b(theta_k)|^p}.
	
	For $E_{2}$ in \eqref{bound dEf(Z)}, we discuss it on set $\{\vert Z_{t}\vert\le R_{1}\}$ first. Using  inequality \eqref{fDoublePrime}, which holds only on $\{\vert Z_{t}\vert\le R_{1}\}$, \eqref{ine Boundsigmay} with noting the concavity of $f$, and \eqref{DissipativeDeduced} to get 
	\begin{align}\label{bound dEf(Z) Pt2}
		\begin{split}
			E_{2}\le&\mathbb{E}\left[\vert Z_{t}\vert\kappa(\vert Z_{t}\vert)f^{\prime}(\vert Z_{t}\vert)+\frac{2}{K_{3}}[\phi_{1}^{\delta}(Z_{t})]^{2}f^{\prime\prime}(\vert Z_{t}\vert)\right]\\
			\le&\mathbb{E}\left[-c_{1}[\phi_{1}^{\delta}(Z_{t})]^{2}f(\vert Z_{t}\vert)+[\phi_{2}^{\delta}(Z_{t})]^{2}\vert Z_{t}\vert\kappa(\vert Z_{t}\vert)f^{\prime}(\vert Z_{t}\vert)-2c_{2}[\phi_{1}^{\delta}(Z_{t})]^{2}\right]\\
			=&\mathbb{E}\left[-c_{1}f(\vert Z_{t}\vert)-2c_{2}+[\phi_{2}^{\delta}(Z_{t})]^{2}\left(c_{1}f(\vert Z_{t}\vert)+2c_{2}\right)+[\phi_{2}^{\delta}(Z_{t})]^{2}\vert Z_{t}\vert\kappa(\vert Z_{t}\vert)f^{\prime}(\vert Z_{t}\vert)\right]\\
			\le& \mathbb{E}\left[-c_{1}f(\vert Z_{t}\vert)-2c_{2}+c_{1}\delta+\sup_{r\in[0,\delta]}r\kappa(r)+2c_{2}1_{\{\vert Z_{t}\vert<\delta\}}\right],
		\end{split}
	\end{align} 
	where the equality follows from \eqref{Phi1+Phi2=1}, and \eqref{ThreeUpperBounds} is used in the last step. Then consider $E_{2}$ on the set $\{\vert Z_{t}\vert> R_{1}\}$. In such case,  we have \eqref{fDoubleprimeGeneral} instead of \eqref{fDoublePrime}, so the following estimation holds instead:
	\begin{align}\label{bound dEf(Z) Pt2 Outside}
		E_{2}\le \mathbb{E}\left[-c_{1}f(\vert Z_{t}\vert)+c_{1}\delta+\sup_{r\in[0,\delta]}r\kappa(r)\right].
	\end{align} 
	Combining \eqref{bound dEf(Z) Pt1}, \eqref{bound dEf(Z) Pt2}, \eqref{bound dEf(Z) Pt2 Outside},  \eqref{bound dEf(Z)}, on $\{\vert Z_{t}\vert\le R_{1}\}$, it follows that\begin{align}\label{dEf(|Z|)}
		\frac{{\rm d}}{{\rm d}t}\mathbb{E}[f(\vert Z_{t}\vert)]\le -C_{1}\mathbb{E} f(\vert Z_{t}\vert)-2c_{2}+c_{1}\delta+\sup_{r\in[0,\delta]}r\kappa(r)+2c_{2}\mathbb{P}\left(\vert Z_{t}\vert<\delta\right)+C_{2}(1+\vert z\vert)\eta_{n}^{1/2}.
	\end{align} 
	 and on $\{\vert Z_{t}\vert> R_{1}\}$, one has\begin{align}\label{dEf(|Z|)Outside}
		\frac{{\rm d}}{{\rm d}t}\mathbb{E}[f(\vert Z_{t}\vert)]\le -C_{1}\mathbb{E} f(\vert Z_{t}\vert)+c_{1}\delta+\sup_{r\in[0,\delta]}r\kappa(r) +C_{2}(1+\vert z\vert)\eta_{n}^{1/2}.
	\end{align}  
		
	Finally, combining \eqref{dE|V(x)+V(y)|} with   $\varepsilon=2q^{-1}c_{2}$,  \eqref{dEf(|Z|)}, and \eqref{dErho(X,Y)}, \eqref{ddtErho1} can be checked on set $\{\vert Z_{t}\vert\le R_{1}\}$. While on set $\{\vert Z_{t}\vert> R_{1}\}$,  we can use \begin{align*}
		{\rm d}\mathbb{E}[\varepsilon V(X_{t})+\varepsilon V(Y_{t})]\le -\frac{\lambda}{2} \mathbb{E}[\varepsilon V(X_{t})+\varepsilon V(Y_{t})]{\rm d}t
	\end{align*}
	instead of \eqref{dE|V(x)+V(y)|}, because $\vert Z_{t}\vert\ge R_{1}$ implies $\vert X_{t}\vert^{2}+\vert Y_{t}\vert^{2}\ge\frac{R_{1}^{2}}{4}$ and one can enlarge $R_{1}$ such that $\varepsilon(8q-\lambda R_{1}^{2})\le 0$ for ensuring the estimation above holds. This, combined with \eqref{dEf(|Z|)Outside}, shows that \eqref{ddtErho1} also holds on $\{\vert Z_{t}\vert\ge R_{1}\}$. 
\end{proof}

\section{Proof of Theorem \ref{TheoremCriterion}}\label{AppendixProofIP}

Here we follow the argument in proving McLeish's FCLT criterion \ocite{McL74}*{Theorem 3.2} to prove its generalization Theorem \ref{TheoremCriterion}.

Begin with truncation. The Lemma \ref{LemmaTruncation} below can be verified in the same way as \ocite{McL74}*{Theorem 3.2}, we omit its proof.
\begin{lemma}\label{LemmaTruncation}
	Under assumptions of Theorem \ref{TheoremCriterion}, let \begin{align}
		\bar{X}_{n,i}\coloneqq X_{n,i}1_{\{\sum_{m=1}^{i-1}X_{n,m}^{2}\le a(T)+1\}}\label{def barXni},
	\end{align}
	then $\{\bar{X}_{n,i}\}_{i=1}^{k_{n}(T)}$ is a martingale difference array satisfying, for each $t\in[0,T]$:
	\begin{align}
		\max_{i\le k_{n}(T)}\vert\bar{X}_{n,i}\vert\to 0\quad\text{ in probability}\label{IP AltCond1 UI}\tag{A1}
	\end{align} 
	and \begin{align}
		\sum_{i=1}^{k_{n}(t)}\bar{X}_{n,i}^{2}\to a(t)\quad \text{ in }L^{1},\quad t\in[0,T].\label{IP AltCond2 L1Conv}\tag{A2}
	\end{align}
	Moreover,
	\begin{align}\label{TruncationIsSafe}
		\mathbb{P}(\exists i\le k_{n}(t),\ \bar{X}_{n,i}\ne X_{n,i})\to 0
	\end{align}
	for each $t\in[0,T]$ as $n\to\infty$.
\end{lemma}
We then show that the truncated array admits finite dimensional distribution convergence and tightness. In particular, we additionally need Assumption \ref{Asp a} for showing tightness, which automatically holds under the setting of  \ocite{McL74}*{Theorem 3.2}. 
\begin{proof}[Proof of Theorem \ref{TheoremCriterion}]
	Consider the truncated martingale difference array $\{\bar{X}_{n,i}\}_{i=1}^{k_{n}(T)}$ as in \eqref{def barXni}, which satisfies \eqref{IP AltCond1 UI} and \eqref{IP AltCond2 L1Conv} as in Lemma \ref{LemmaTruncation}. Denoting 
	\begin{align*}
		\bar{S}_{n}(t)\coloneqq\sum_{i=1}^{k_{n}(t)}\bar{X}_{n,i},
	\end{align*}
	taking $t=T$ in \eqref{TruncationIsSafe} gives that $\mathbb{P}(S_{n}\ne \bar{S}_{n})\to 0$, so it is sufficient to show that $\bar{S}_{n}(\cdot)\Rightarrow B_{a(\cdot)}^{1}$ on $D[0,T]$. To this end, we prove the convergence of finite dimensional distribution and tightness respectively. 
	
	For the convergence of finite dimensional distribution, we use Cram\'{e}r--Wold argument established in \ocite{CW36}. Specifically, consider  arbitrary points $0=t^{0}<t^{1}<\cdots<t^{m}<T$ in $[0,T]$, let $u_{1},\cdots,u_{m}$ be arbitrary reals, put \begin{align*}
		Y_{n,i}\coloneqq\left\{\begin{array}{ll}
			u_{j}\bar{X}_{n,i},&k_{n}(t^{m})\ge i,\\ 
			0,&k_{n}(t^{m})<i,
		\end{array}\right.\quad\text{ where } j\coloneqq\inf\{\ell:k_{n}(t^{\ell})\ge i\}.
	\end{align*}
	Note that \begin{align*}
		\sum_{i=1}^{k_{n}(T)}Y_{n,i}=\sum_{j=1}^{m}u_{j}[\bar{S}_{n}(t^{j})-\bar{S}_{n}(t^{j-1})],\quad \sum_{i=1}^{k_{n}(T)}Y_{n,i}^{2}=\sum_{j=1}^{m}u_{j}^{2}\left[\sum_{i=1}^{k_{n}(t^{j})}\bar{X}_{n,i}^{2}-\sum_{i=1}^{k_{n}(t^{j-1})}\bar{X}_{n,i}^{2}\right].
	\end{align*}
	Denote $\sigma_{u}^{2}\coloneqq\sum_{j=1}^{m}u_{j}^{2}(a(t^{j})-a(t^{j-1}))$. Condition \eqref{IP AltCond2 L1Conv} gives that  \begin{align}\label{CorCond1}
		\sum_{i=1}^{k_{n}(T)}\mathbb{E} Y_{n,i}^{2}=\sum_{j=1}^{m}u_{j}^{2}\left[\mathbb{E}\sum_{i=1}^{k_{n}(t^{j})}\bar{X}_{n,i}^{2}-\mathbb{E}\sum_{i=1}^{k_{n}(t^{j-1})}\bar{X}_{n,i}^{2}\right]\to \sigma_{u}^{2}
	\end{align}
	as $n\to\infty$. On the other hand, by Markov's inequality and $\eqref{IP AltCond2 L1Conv}$, one has \begin{align}
		\label{CorCond2}\mathbb{P}\left(\sigma_{u}^{2}-\sum_{i=1}^{k_{n}(T)}Y_{n,i}^{2}>\varepsilon\right)\le\mathbb{P}\left(\left\vert\sigma_{u}^{2}-\sum_{i=1}^{k_{n}(T)}Y_{n,i}^{2}\right\vert>\varepsilon\right)\le \varepsilon^{-1}\mathbb{E}\left\vert\sigma_{u}^{2}-\sum_{i=1}^{k_{n}(T)}Y_{n,i}^{2}\right\vert \to0
	\end{align}
	as $n\to\infty$. Under \eqref{CorCond1} and \eqref{CorCond2},  \ocite{McL74}*{Corollary 2.8}  applies, which gives the following weak convergence as $n\to\infty$:
	\begin{align*}
		\sum_{j=1}^{m}u_{j}[\bar{S}_{n}(t^{j})-\bar{S}_{n}(t^{j-1})]=\sum_{i=1}^{k_{n}(T)}Y_{n,i} \Rightarrow\mathcal{N}(0,\sigma_{u}^{2})\stackrel{d}{=}\sum_{i=1}^{m}u_{i}[B^{1}(a(t^{i}))-B^{1}(a(t^{i-1}))].
	\end{align*}
	Convergence of the finite dimensional distribution been established.
	
	To show the tightness, as in showing  \ocite{McL74}*{Theorem 3.2} , it is sufficient to show for each $\varepsilon>0$ that \begin{align*}
		\lim_{\delta\to 0}\varlimsup_{n\to\infty}\mathbb{P}\left(\sup_{\substack{\vert s-t\vert\le \delta,\\s,t\in[0,T]}}\vert\bar{S}_{n}(s)-\bar{S}_{n}(t)\vert>\varepsilon\right)=0.
	\end{align*}
	To derive this limit, we follow the approach in proving  \ocite{Bro71}*{Theorem 3} . Recall notations in Assumption \ref{Asp a}. By considering a $\delta$-mesh on interval $[0,T]$, inserting terms $\bar{S}_{n}(t_{\delta}^{m})$ and taking supremum over $m$, it follows that \begin{align*}
		\sup_{\substack{\vert s-t\vert\le \delta,\\s,t\in[0,T]}}\vert\bar{S}_{n}(s)-\bar{S}_{n}(t)\vert  \le4\sup_{0\le m\le \lfloor \delta^{-1}T\rfloor}\sup_{t\in[t_{\delta}^{m},t_{\delta}^{m+1}]}\vert \bar{S}_{n}(t)-\bar{S}_{n}(t_{\delta}^{m})\vert.
	\end{align*}
	Therefore, one has \begin{align*}
		\mathbb{P}\left(\sup_{\substack{\vert s-t\vert\le \delta,\\s,t\in[0,T]}}\vert\bar{S}_{n}(s)-\bar{S}_{n}(t)\vert>\varepsilon\right)\le&\sum_{m=0}^{\lfloor \delta^{-1}T\rfloor}\mathbb{P}\left(\sup_{t\in[t_{\delta}^{m},t_{\delta}^{m+1}]}\vert \bar{S}_{n}(t)-\bar{S}_{n}(t_{\delta}^{m})\vert>\frac{\varepsilon}{4}\right)\\
		=&\sum_{m=0}^{\lfloor \delta^{-1}T\rfloor}\mathbb{P}\left(\sup_{t\in[t_{\delta}^{m},t_{\delta}^{m+1}]}\left\vert\sum_{i=k_{n}(t_{\delta}^{m})}^{k_{n}(t)}\bar{X}_{n,i}\right\vert>\frac{\varepsilon}{4}\right)\\
		\le&\frac{8}{\varepsilon}\sum_{m=0}^{\lfloor \delta^{-1}T\rfloor}\mathbb{E}\left[\vert Z^{(n)}_{m}\vert 1_{\{\vert Z^{(n)}_{m}\vert> \varepsilon/8\}}\right],
	\end{align*}
	where $Z_{m}^{(n)}\coloneqq\bar{S}_{n}(t_{\delta}^{m+1})-\bar{S}_{n}(t_{\delta}^{m})$ and the last step follows from a maximal inequality for martingales; see \ocite{Bro71}*{Lemma 4}. 
	
	The convergence of finite dimensional distribution gives $Z_{m}^{(n)}\Rightarrow Z_{m}$ as $n\to\infty$, where $Z_{m}\stackrel{d}{=}\mathcal{N}(0,\tau_{m}(\delta))$ and $\tau_{m}(\delta)$  is defined as in   Assumption \ref{Asp a}. By Skorohod representation theorem \ocite{Bil79}*{Theorem 6.7}, one may find $\{\widetilde{Z}_{m}^{(n)}\}_{n\ge 1}$, $\widetilde{Z}_{m}$ with $\widetilde{Z}_{m}^{(n)}\stackrel{d}{=}Z_{m}^{{(n)}}$ and  $\widetilde{Z}\stackrel{d}{=}Z_{m}$ on some probability space, such that $\widetilde{Z}_{m}^{(n)}\to Z_{m}$ a.s. as $n\to\infty$. Note by condition \eqref{IP AltCond2 L1Conv} that \begin{align*}
		\mathbb{E}\vert \widetilde{Z}_{m}^{(n)}\vert^{2}=\mathbb{E}\vert Z_{m}^{(n)}\vert^{2}=\mathbb{E}\sum_{i=1}^{k_{n}((m+1)\delta)}\bar{X}_{n,i}^{2}-\mathbb{E}\sum_{i=1}^{k_{n}(m\delta)}\bar{X}_{n,i}^{2}\to \tau_{m}(\delta)=\mathbb{E}\vert Z_{m}\vert^{2}<\infty, 
	\end{align*} 
	which gives uniformly integrability of $\{\widetilde{Z}_{m}^{(n)}\}_{n\ge 1}$. By using the Vitali convergence theorem \ocite{Sch17}*{Theorem 16.6}, one has $\widetilde{Z}_{m}^{(n)}\to \widetilde{Z}_{m}$ in $L^{1}$ as $n\to\infty$, which implies \begin{align*}
		\lim_{n\to\infty}\mathbb{E}[\vert Z^{(n)}_{m}\vert 1_{\{\vert Z^{(n)}_{m}\vert> \varepsilon/8\}}]=\mathbb{E}[\vert Z_{m}\vert 1_{\{\vert Z_{m}\vert> \varepsilon/8\}}].
	\end{align*}
	Therefore, one arrives at \begin{align*}
		\varlimsup_{n\to\infty}	\mathbb{P}\left(\sup_{\substack{\vert s-t\vert\le \delta,\\s,t\in[0,T]}}\vert\bar{S}_{n}(s)-\bar{S}_{n}(t)\vert>\varepsilon\right)\le&\frac{C}{\varepsilon}\sum_{m=0}^{\lfloor \delta^{-1}T\rfloor}\int_{\{\vert Z_{m}\vert>\varepsilon/8\}}\vert Z_{m}\vert {\rm d}\mathbb{P}\\
		\le&\frac{C}{\varepsilon^{2}}\sum_{m=0}^{\lfloor \delta^{-1}T\rfloor}\int_{\{\vert Z_{m}\vert>\varepsilon/8\}}\vert Z_{m}\vert^{2} {\rm d}\mathbb{P}\\
		=&\frac{C}{\varepsilon^{2}}\sum_{m=0}^{\lfloor \delta^{-1}T\rfloor}\sqrt{\frac{\tau_{m}(\delta)}{2\pi}}\int_{\{\vert x\vert>\varepsilon/8\sqrt{\tau_{m}(\delta)}\}}x^{2}e^{-\frac{x^{2}}{2}} {\rm d} x\\
		\le&C\delta^{-1}\sqrt{\tau(\delta)}\exp\left(-C \tau(\delta)^{-1}\right),
	\end{align*}
	which tends to $0$ as $\delta\to0$ by Assumption \ref{Asp a}. The tightness follows.
\end{proof}

	\bibliography{./ref}
\end{document}